\documentclass{article}

\usepackage{PRIMEarxiv}

\usepackage{lineno,hyperref}
\usepackage{compactbib}
\usepackage{amsmath}
\usepackage{amsfonts}
\usepackage{amssymb}
\usepackage{txfonts}
\usepackage{graphics}
\usepackage{graphicx}
\usepackage{epsfig}
\usepackage{color}
\usepackage{xcolor}
\usepackage{fleqn}
\usepackage{float}
\usepackage{bm}

\usepackage{algorithm}
\PassOptionsToPackage{noend}{algpseudocode}
\usepackage{algpseudocode}

\errorcontextlines\maxdimen

\usepackage{PRIMEarxiv}

\usepackage[utf8]{inputenc} 
\usepackage[T1]{fontenc}    
\usepackage{hyperref}       
\usepackage{url}            
\usepackage{booktabs}       
\usepackage{nicefrac}       
\usepackage{microtype}      
\usepackage{fancyhdr}       

\title{A change of measure enhanced near exact Euler Maruyama scheme for the solution to nonlinear stochastic dynamical systems}

\author{Tapas Tripura\\
  Department of Applied Mechanics\\
  Indian Institute of Technology Delhi\\
  \texttt{tapas.t@am.iitd.ac.in} \\

   \And
  Mohammad Imran \\
  Department of Civil Engineering\\
  Indian Institute of Technology Guwahati\\
  \texttt{mimran@iitg.ac.in}
  
  \And
  Budhaditya Hazra \\
  Department of Civil Engineering\\
  Indian Institute of Technology Guwahati\\
  \texttt{budhaditya.hazra@iitg.ac.in} \\
  
  \And
      Souvik Chakraborty \\
  Department of Applied Mechanics\\
  School of Artificial Intelligence (ScAI)\\
  Indian Institute of Technology Delhi\\
  \texttt{souvik@am.iitd.ac.in} \\
}

\begin{document}
\maketitle

\begin{abstract}
The present study utilizes the Girsanov transformation based framework for solving a nonlinear stochastic dynamical system in an efficient way in comparison to other available approximate methods. In this approach, a rejection sampling is formulated to evaluate the Radon-Nikodym derivative arising from the change of measure due to Girsanov transformation. The rejection sampling is applied on the Euler Maruyama approximated sample paths which draw exact paths independent of the diffusion dynamics of the underlying dynamical system. The efficacy of the proposed framework is ensured using more accurate numerical as well as exact nonlinear methods. Finally, nonlinear applied test problems are considered to confirm the theoretical results. The test problems demonstrates that the proposed exact formulation of the Euler-Maruyama provides an almost exact approximation to both the displacement and velocity states of a second order non-linear dynamical system.
\end{abstract}


\section{Introduction}
Modeling engineering dynamical systems and determining their response is becoming more versatile as new approaches for solving nonlinear engineering systems are being discovered frequently. Since engineering models are typically uncertain in nature, modeling errors or naturally occurring noise terms are inevitable. As a result, research on numerical approaches for the approximation of stochastically driven oscillators has been a focused area \cite{gikhman2007stochastic,oksendal2013stochastic}. Both linear and non-linear oscillators related to physics and engineering fields are reducible to first order stochastic differential equations (SDEs) through suitable state transformation. In general, SDEs are extremely difficult to solve analytically. Although for some cases it is possible to find stationary density functions of the oscillators by using Fokker Planck Equation (FPK) \cite{maruyama1955continuous}, many of the approaches for determining the analytical solution of stochastically driven nonlinear oscillators suffer from the curse of dimensionality \cite{lin1988exact}. Since analytical steady-state solutions are not always possible, therefore their existence imposes limits on exploratory dynamics of the oscillators. As a result, many attempts have been made to develop efficient numerical methods. Monte Carlo simulation (MCS) aided with direct stochastic numerical integration is a more extensively used but less elegant for solving issues of substantially higher complexity as the dimension of the problem increases. Nonetheless, the precision of direct integration schemes and the highly repetitive computation over a potentially large ensemble limit the utility of MCS approaches. Numerical integration techniques for stochastically driven oscillators are frequently constructed using stochastic Taylor expansions, which have a computational disadvantage due to the difficulties of computing numerous multiple stochastic integrals (MSIs) \cite{kloeden1992stochastic}. A few schemes in this category are Euler-Maruyama \cite{maruyama1955continuous}, Milstein \cite{milstein1994numerical}, stochastic Heun \cite{gard1988introduction}, stochastic Runge-Kutta \cite{ruemelin1982numerical}, Strong Taylor 1.5 \cite{kloeden1992higher} and Weak order 3.0 Taylor $O(\Delta t^{3})$ \cite{tripura2020ito} etc. Being the simplest and easiest among the available stochastic integration techniques, Euler-Maruyama (EM) is of primary interest and an improvement based on change of measure is the aim of this work.

There exist approximate analytical methods which can determine the statistical responses of the stochastic systems. For reference Socha \cite{socha1999statistical} has provided a comprehensive review of equivalent and statistical linearization of nonlinear oscillators. Generally, most linearization methods such as Equivalent Linearization, higher order linearization \cite{iyengar1988higher}, equivalent non-linearization \cite{caughey1986response}, stochastic averaging \cite{liu2012stochastic}, Gaussian equivalent linearization\cite{anh2003improved} globally replace the nonlinear drift field by an equivalent, time-variant and linear function such that the error in the replacement is minimized in a mean square sense. All of these methods use iterative approaches, however, it should be emphasized that the scope of global linearization for nonlinear oscillators are relatively limited and can become unacceptably imprecise at times since obtaining the linearized coefficients of the linear function is dependent on the order of optimization \cite{burt1981quantum}.

Another numerical approach for solving non-linear systems under additive stochastic excitation is the Girsanov based formulation of the local transversal linearization (LTL) schemes \cite{biscay1996local,raveendran2013nearly}. In this family of methods, the errors in the non-linearity approximation by local linearization schemes are absorbed inside the stochastic diffusion term. The modification in the diffusion is evaluated using the change of measure of the process referred to as the Girsanov transformation. The change of measure arising due to this adjustment in diffusion helps in formulating a correction process that improves the low order of convergence of LTL schemes. The change of measure applied over probability measures is related by a density function known as the Radon-Nikodym derivative. This Radon-Nikodym derivative can be evaluated directly through It\^{o}-Taylor expansion of the stochastic exponential or through stochastic filtering and rejection sampling. These schemes have shown efficacy and early promise over It\^{o}-Taylor based numerical integration schemes, however, for higher dimensional problem yet to be examined. Inspite of the advantage over higher order It\^{o}-Taylor schemes, the LTL schemes are computationally expensive due to evaluation of stochastic matrix exponential. The present work emphasizes development of more efficient and accurate yet simple alternatives utilizing the change of measure. 

In the case of EM method, simulating a multi-dimensional stochastic system using the EM discretization is a simple and efficient numerical task since it does not involve computation of any matrix exponential. Apart from the computational inefficiency, the stochastic matrix exponential involve stochastic Brownian integrals which are not straightforward to evaluate and EM scheme is advantageous in such situations. Many efficient improvements in the EM method have been done for solving different types of SDEs such as composite EM \cite{burrage2001composite}, implicit Euler–Taylor \cite{tian2001implicit}, truncated EM method \cite{mao2015truncated}, composite previous-current-step EM \cite{nouri2018improved}. Inspite of the modifications, these schemes are still approximate in nature. A new framework for near exact simulation of SDEs using EM is presented in this study which overcomes the absence of higher order terms from It\^{o}-Taylor expansion in EM approximation and even for high dimensional systems assumes no approximation and provides highly efficient computation. In the proposed framework, a Girsanov's change of measure is adopted over the EM discretization, where the Radon-Nikodym derivative is evaluated using rejection sampling \cite{beskos2005exact}. The change of measure between two probability measures generates a dynamic probability density function at every ensemble step and when the rejection sampling is performed appropriately, it returns near exact draws from any finite-dimensional distribution of the EM solution for the SDE in interest. Since, the mathematical aspect of the process dynamics is handled by the EM, the draws from rejection sampling is independent of the diffusion dynamics. Altogether, the idea then centers around the restriction that the solution paths of EM at each step must satisfy the density function (or, Radon-Nikodym derivative). Finally, one can achieve both computationally efficient and near exact solutions through the proposed framework, since sample generation using EM requires less computational time and the rejection sampling on the generated samples can be done suitably by defining appropriate events following standard approaches \cite{beskos2005exact}. 

From this onward the paper is arranged as follows: \textbf{Section 2}: a short background on types of SDE and change of measure, the way it serves the current interest is provided. \textbf{Section 3}: the proposed change of measure framework  for Euler-Maruyama and rejection sampling algorithm for evaluating Radon-Nikodym derivative is presented. \textbf{Section 4}: numerical illustration using a fairly representative class of non-linear oscillators excited using zero mean Gaussian white noise are provided. \textbf{Section 5}: the paper is concluded by highlighting key achievements of the work. 

\section{Background on SDE and Radon-Nikodym derivative for Brownian process}

\subsection{Euler-Maruyama approximation of Stochastic differential equation}\label{EM}
Let the probability space $\left( {\Omega ,\mathcal{F},P} \right)$, with natural filtration $({{({\mathcal{F}_t})},0 \le t \le T)}$ be constructed from sub $\sigma$-algebras of $\mathcal{F}$. Consider, an $m$-dimensional $n$-factor SDE which has a deterministic dynamics modeled as drift driven by the additive volatile component:  
\begin{equation}\label{sdeg}
d{\bm{X}}_t = {{\bf g}}\left( {t,{{\bm{X}}_t}} \right)dt + \sum\limits_{j = 1}^n {{\bf f}_j\left( {t,{{\bm{X}}_t}} \right)} d{{\bm W}_j}\left( t \right); \qquad
{\bm X}(t=t_0)={\bm X}_0; \qquad t \in [0,T]
\end{equation}
where, ${{\bm{X}}_t} \in {\mathbb{R}^m}$ denotes the ${{\mathcal{F}_t}}$-measurable state vector, ${\bf{g}}\left( {t,{{\bm{X}}_t}} \right) \in \left[ {0,T} \right] \times {\mathbb{R}^m} \mapsto {\mathbb{R}^m}$ is drift function and ${{\bf{f}}_j}\left( {t,{{\bm{X}}_t}} \right) \in \left[ {0,T} \right] \times {\mathbb{R}^m} \mapsto {\mathbb{R}^m}$ for $j = 1, \ldots n$ is volatility coefficient function. Here, $\{{\bm W}_j(t):j = 1, \ldots n\}$ is the independent Brownian process with respect to the probability measure $P$. The bound and uniqueness of the solution vector ${{\bf{X}}_t}$ are defined by the  following criterion \cite{hassler2016stochastic}:
\begin{enumerate}
	\item Lipschitz continuity condition: the functions ${\bf g}\left( {t,{{\bm {X}}_t}} \right)$ and ${\bf f}\left( {t,{{\bm {X}}_t}} \right)$ must be partially differentiable with respect to ${{\bm{X}}_t}$,
	\begin{equation}
	\left| {{\bf{g}}\left( {t,{{\bm{X}}_t}} \right)} - {{\bf{g}}\left( {t,{{\bm{Y}}_t}} \right)} \right| + \left| {{\bf{f}}\left( {t,{{\bm{X}}_t}} \right)} - {{\bf{f}}\left( {t,{{\bm{Y}}_t}} \right)} \right| \le {\mathcal{D}}\left| {\bm{X}} - {\bm{Y}} \right|
	\end{equation}
	\item Boundness condition on the growth of the diffusion process: 
	\begin{equation}
	\left| {{\bf{g}}\left( {t,{{\bm{X}}_t}} \right)} \right| + \left| {{\bf{f}}\left( {t,{{\bm{X}}_t}} \right)} \right| \le {\mathcal{C}_1} + {\mathcal{C}_2}\left| {\bm{X}} \right|
	\end{equation}
	\item For a well defined initial point ${{\bm{X}}({t_0})}$:
	\begin{equation}
	E\left[ {{{\left| {{\bm{X}}({t_0})} \right|}^2}} \right] < \infty 
	\end{equation}
\end{enumerate}
where, ${{\bm{X}}(t)}$ and ${{\bm{Y}}(t)}$ are real variables and ${\mathcal{D}}$, ${\mathcal{C}_1}$ and ${\mathcal{C}_2}$ are some positive constants. Under the uniqueness property of the solution, the Euler-Maruyama (EM) approximation to the time evolution of the SDE in Eq. (\ref{sdeg}) over the time interval $t \in \left[ {{t_{i - 1}},{t_i}} \right]$, is given as \cite{maruyama1955continuous}:
\begin{equation}\label{explicit}
{{\bm{X}}_{{t_i}}} = {{\bm{X}}_{{t_{i - 1}}}} + {\bf{g}}\left( {{t_{i - 1}},{{\bm{X}}_{{t_{i - 1}}}}} \right)\Delta t + \sum\limits_{j = 1}^n {{\bf{f}}\left( {{t_{i - 1}},{{\bm{X}}_{{t_{i - 1}}}}} \right)} \Delta {{\bm{W}}_j}\left( {{t_i}} \right)
\end{equation}
where, ${{\bm{X}}_{{t}}}$ represents the EM based approximated state vector, $\Delta t$=$(t_i - t_{i-1})$ is the time increment and $\Delta {{\bm{W}}_{{t_i}}} = \left( {{{\bm{W}}_{{t_i}}} - {{\bm{W}}_{{t_{i - 1}}}}} \right)$ is the Brownian increment having a Gaussian distribution $\sim \mathcal{N}\left( {0,\Delta t} \right)$. In the absence of multiple Brownian integrals of higher order of smallness, the EM scheme has a strong order of convergence of $O ({\Delta t}^{0.5})$ and a weak order $O ({\Delta t}^{1})$, which can be verified from the order of $\Delta {\bm{W}}$ and $\Delta {t}$, respectively.

\subsection{Change of measure}
Let $X \sim \mathcal{N}[0,1]$ be a standard normal variable under the probability space be $\left( {\Omega ,\mathcal{F},P} \right)$, whose density and distribution functions are represented by $p(x)$ and $P(x)$. The density and distribution functions are related by the relation: $dP(x) = p(x)dx$. Let us define an event A:= ${\alpha  \le X \le \beta }$ under the probability measure $P$. Further let us assume that a shift in the intervals by a constant $\mu$ is applied, which is formally represented as:
\begin{equation}
P[\underbrace {\alpha  - \gamma  \le X \le \beta  - \gamma }_{A - \gamma }] = \frac{1}{{\sqrt {2\pi } }}\int_{\alpha  - \gamma }^{\beta  - \gamma } {\exp \left( { - \frac{1}{2}{x^2}} \right)dx} 
\end{equation}
Without invoking the effect of the shift in the interval to the law of evolution of random variable $X$, a new probability measure can suitably be determined through modification of the distribution function. A simple mathematical substitution of variables yield:
\begin{equation}
P[\alpha  - \gamma  \le X \le \beta  - \gamma ] = \frac{1}{{\sqrt {2\pi } }}\int_\alpha ^\beta  {\exp \left( { - \frac{1}{2}{{(x - \gamma )}^2}} \right)dx}  = \int_\alpha ^\beta  {\exp\left( {x\gamma  - \frac{1}{2}{\gamma ^2}} \right)dP(x)} 
\end{equation}
In the above equation, it can be seen that the shift in the intervals is absorbed into the probability distribution function. The shift in the distribution function causes the change of measure $P \to Q$ as:
\begin{equation}
P\left[ {A - \gamma } \right] = \int\limits_A {\exp \left( {x\gamma  - \frac{1}{2}{\gamma ^2}} \right)dP(x)}  = {E^p}\left[ {\exp \left( {x\gamma  - \frac{1}{2}{\gamma ^2}} \right){I_{[\alpha ,\beta ]}}} \right] = {E^p}\left[ {\Lambda (x){I_A}} \right] = Q[A]
\end{equation}
where, ${E^p}[.]$ is the expectation with respect to the probability measure P and ${{I_A}}$ is the indicator function which returns 1 if X lies in the specified interval [${\alpha ,\beta}$], otherwise 0. The above equation means, the RV $X - \gamma $ under probability measure $Q \sim \mathcal{N}\left[ {0,1} \right]$ has identical property as in measure $P$. The term $\Lambda (x) = \exp \left( {x\gamma  - \frac{1}{2}{\gamma ^2}} \right)$ defines the relationship between two probability measures $P$ and $Q$, often called as Radon-Nikodym derivative:
\begin{equation}
\frac{{dQ(x)}}{{dP(x)}} = \Lambda (x)
\end{equation}
The Radon-Nikodym derivative is indeed a density function with the property $\Lambda (x) > 0$ and ${E^P}\left[ {\Lambda (x)} \right] = 1$. The former one can be easily verified since $\Lambda (x)$ is a exponential function of $X$ and the later one is proved by taking expectation of $\Lambda (x)$ in the measure $P$. To calculate the probability of the path of the random variable (RV) $X$, let us consider the path ${\bm{X}} = \left\{ {X(t):0 \le t \le {t_N},X(t_0) = 0} \right\}$ for $N$-partitions between the interval $t \in \left[ {0,{t_N}} \right]$. With the initial condition ${X}\left( {{t_0}} \right) = 0$, let us assume that the increments of the RV $X$ are independent and identically distributed (i.i.d) with normal distribution having mean $0$ and variance $\Delta t$ i.e. $\Delta {X_i} \sim {\cal N}[X({t_0}),\Delta t]$. If the incremental process is given as: ${X_i} = {X_{i - 1}} + \Delta {X_i}$, then, the increments will have normal distribution with $\{\Delta {X_1} \sim N[{X_0},\Delta {t_1}],\Delta {X_2} \sim N[{X_1},\Delta {t_2}], \ldots ,\Delta {X_N} \sim N[{X_{N - 1}},\Delta {t_N}]\}$. For uniform $\Delta t$ this can be represented by $\Delta {X_i} \sim N[{X_{i-1}},\Delta t]$. 
With this knowledge, given a starting value the density function of the process at a target point $X_i$ over the interval $\Delta t$ can be expressed as a conditional distribution as follows,
\begin{equation}\label{density}
P({X_i} - {X_{i - 1}}\left| {{X_{i - 1}}} \right.) = \frac{1}{{\sqrt {2\pi \Delta t} }}\exp \left( { - \frac{1}{{2\Delta t}}{{({x_i} - {x_{i - 1}})}^2}} \right)
\end{equation}
In order to find the density function of the complete path of the RV $X$, one simply needs to find the $N$-product of the density function in Eq. (\ref{density}), since the increments are i.i.d. Then the complete probability of the process can be found by performing integration over the density function. Towards this, let us consider an event A as, $A := \left[ {\left( {{\alpha _1} \le {X_1} \le {\beta _1}} \right), \ldots ,\left( {{\alpha _N} \le {X_N} \le {\beta _N}} \right)} \right]$. For large $N$, the distribution for the event A can be represented by following infinite dimensional integral over the $N$-product of density functions in Eq. (\ref{density}):
\begin{equation}
P[A] = \frac{1}{{\sqrt {2\pi \Delta t} }}\int_{{\alpha _1}}^{{\beta _1}} {\int_{{\alpha _2}}^{{\beta _2}} { \ldots \int_{{\alpha _N}}^{{\beta _N}} {\exp\left( { - \frac{1}{{2\Delta t}}\sum\limits_{i = 1}^N {{{({x_i} - {x_{i - 1}})}^2}} } \right)d{x_1}d{x_2} \ldots d{x_N}} } } 
\end{equation}
A stochastic process for the RV $X$ under the above finitely large collection of densities is guaranteed due to the Kolmogorov extension theorem \cite{oksendal2013stochastic}. To introduce the change of measure, let us apply a shift in the event A as follows:
\begin{equation}
\left[ {A - {y_k}} \right] = \left[ {\left( {\left( {{\alpha _1} - {\mu _1}\Delta t} \right) \le {X_1} \le \left( {{\beta _1} - {\mu _1}\Delta t} \right)} \right), \ldots ,\left( {\left( {{\alpha _N} - \sum {{\mu _i}\Delta t} } \right) \le {X_N} \le \left( {{\beta _N} - \sum {{\mu _i}\Delta t} } \right)} \right)} \right]
\end{equation}
The shift, ${y_k} = \sum_{i = 1}^N {{\gamma _i}\Delta t} $, in the event A introduces mutation in the distribution function, which can be expressed in the limits of original event A as the following:
\begin{equation}
\centering
P[A - {y_k}] = \frac{1}{{\sqrt {2\pi \Delta t} }}\int\limits_{{\alpha _1}}^{{\beta _1}} {\int\limits_{{\alpha _2}}^{{\beta _2}} {\ldots \int\limits_{{\alpha _N}}^{{\beta _N}} {\exp\left( { - \frac{1}{{2\Delta t}}\sum\limits_{i = 1}^2 {{{(\Delta {x_i})}^N}} } \right)\underbrace {\exp\left( {\sum\limits_{i = 1}^N {{\gamma _i}\Delta {x_i}}  - \frac{1}{{2\Delta t}}\sum\limits_{i = 1}^N {{\gamma _i}^2\Delta t} } \right)}_{\Lambda ({t_N})}} d{x_1}d{x_2} \ldots d{x_N}} } 
\end{equation}
Noting, $P[A - {y_k}] = \int_A {{\Lambda_n}dP(x)} $ and defining a new measure $Q[A]=P[A - {y_k}]$, the Radon-Nikodym derivative for a multi-dimensional random variable can be defined as:
\begin{equation}\label{2.19}
{\Lambda _N} = \exp\left( {\sum\limits_{i = 1}^N {{\gamma _i}\Delta {x_i}}  - \frac{1}{{2\Delta t}}\sum\limits_{i = 1}^N {{\gamma _i}^2\Delta t} } \right)
\end{equation}
As the number of partitions N$\to \infty$ in the interval $\left[ {{t_0} \ldots {t_N}} \right]$, the Radon Nikodym derivative in the Eq. (\ref{2.19}) becomes a integral:
\begin{equation}\label{2.20}
\sum\limits_{i = 1}^N {{\gamma _i}\Delta {x_i}}  - \frac{1}{{2\Delta t}}\sum\limits_{i = 1}^N {{\gamma _i}^2\Delta t}  \equiv \int_0^t {{\gamma _s}d{x_s} - \frac{1}{2}\int_0^t {{\gamma _s}^2ds} } 
\end{equation}
Thus, the measure $P$ process ${X_N}$ under probability measure $Q$ tends to a distribution $\sim \mathcal{N}\left[ {\sum_{i = 1}^N {{\mu _i}\Delta t} ,\Delta t} \right]$. This also refers to the drift removed process ${X_N} - \sum_{i = 1}^N {{\mu _i}\Delta t} $ under the measure $P$ with the distribution $\sim \mathcal{N}\left[ {0,\Delta t} \right]$ under the probability measure $Q$.

\subsection{Application of change of measure to SDEs}
Consider that under the probability space $\left( {\Omega ,\mathcal{F},P} \right)$, there exists an SDE of the form Eq. (\ref{sdeg}) and ${{\bm{B}}_t} = \left( {\Omega ,\mathcal{F},{{({\bm B}_t)}_t},P} \right)$ be an $n$-dimensional Brownian process under the probability measure $P$. It is to be noted that a shift in the drift under the measure $P$, when applied to the SDE in Eq. (\ref{sdeg}) results an equivalent SDE in a new measure $Q$. This can be formalized by absorbing the shift in the diffusion term as follows:
\begin{equation}\label{shift_brownian}
d{{\bm{X}}_t} = {\bf{g}}\left( {t,{{\bm{X}}_t}} \right)dt + {\bf{f}}\left( {t,{{\bm{X}}_t}} \right)d{\bm{B}}\left( t \right) \to {\bf{f}}\left( {t,{{\bm{X}}_t}} \right)\left( {\frac{{{\bf{g}}\left( {t,{{\bm{X}}_t}} \right)}}{{{\bf{f}}\left( {t,{{\bm{X}}_t}} \right)}}dt + d{{\bm{B}}_t}} \right) = {\bf{f}}\left( {t,{{\bm{X}}_t}} \right)d{{{\tilde{\bm B}}}_t}
\end{equation}
Let, ${{\bf \gamma} \left( {t,{{\bm{X}}_t}} \right)}$ is a $\mathbb{R}^m$-valued process as defined before and there exists a progressively measurable $m$-dimensional process ${{\bf{\gamma }}\left( {s,{{\bm{X}}_s}} \right)}$ with the density $\Lambda_t$ on measure $Q$ under the probability measure $Q$, ($dQ = {{\bf Z} _t}dP$). Then the associated Radon Nikodym derivative for this change of measure $P \to Q$ is identified by substituting ${\bm X}(t)={\bm B}(t)$ as  (\cite{oksendal2013stochastic,baldi2017introduction}):
\begin{equation}\label{radon}
{\bf{\Lambda }}(t) = \frac{{dQ}}{{dP}} = \exp \left( {\int_{{t_{i - 1}}}^t {{\bf{\gamma }}\left( {s,{{\bm{X}}_s}} \right)} d{B_s} - \frac{1}{2}\int_{{t_{i - 1}}}^t {{{\left| {{\bf{\gamma }}\left( {s,{{\bm{X}}_s}} \right)} \right|}^2}} ds} \right)
\end{equation}
In Eq. (\ref{shift_brownian}), ${{\tilde {\bm B}}_t}$ is the Brownian motion in measure $Q$, whose evolution is given by the differential equation: $d{{{\tilde{\bm B}}}_t} = {\bf{\gamma }}\left( {t,{{\bm{X}}_t}} \right)dt + d{{\bm{B}}_t}$, where ${\bf{\gamma }}\left( {t,{{\bm{X}}_t}} \right) = \frac{{{\bf{g}}\left( {t,{{\bm{X}}_t}} \right)}}{{{\bf{f}}\left( {t,{{\bm{X}}_t}} \right)}}$. Here, it can be noticed that the change of measure is purely due to the shift: ${\bf{\gamma }}\left( {t,{{\bm{X}}_t}} \right)dt$. Effected by the shift, the evolution of $Q$-Brownian motion ${{\tilde {\bm B}}_t}$ follows the following rule:
\begin{equation}\label{qbrownian}
{{\tilde {\bm B}}_t} = {{\bm B}_t} + \int_{{t_{i - 1}}}^{{t_i}} {\gamma \left( {s,{{\bm{X}}_s}} \right)ds}; \quad {\tilde{\bm B}}(t=t_0)=0; \quad t\ge0
\end{equation}
It is easy to check that by substituting the relation $d{{\bm{B}}_t} = d{{{\tilde{\bm B}}}_t} - {\bf{\gamma }}\left( {t,{{\bm{X}}_t}} \right)dt$ in Eq. (\ref{sdeg}) will also yield the SDE in $Q$-measure. As a result the SDE becomes,
\begin{equation}
d{{\bm{X}}_t} = {\bf{g}}\left( {t,{{\bm{X}}_t}} \right)dt + {\bf{f}}\left( {t,{{\bm{X}}_t}} \right)\left( {d{{{\tilde{\bm B}}}_t} - \frac{{{\bf{g}}\left( {t,{{\bm{X}}_t}} \right)}}{{{\bf{f}}\left( {t,{{\bm{X}}_t}} \right)}}dt} \right) = {\bf{f}}\left( {t,{{\bm{X}}_t}} \right)d{{{\tilde{\bm B}}}_t}
\end{equation}
The above discussion provides an idea of change of measure theory for Brownian motion, which can be suitably manipulated in a new probability measure  to modify the drift for ease in treatment. This paves way to the idea of how the EM approximation error for the non-linear terms in the drift can be absorbed in the diffusion and correspondingly a change of measure can be formulated for the modified Brownian motion. 

\section{A change of measure framework for near Exact Euler Maruyama (n-EEM)}
In the aforementioned sections, it is understood that the convergence of Euler Maruyama (EM) scheme in its classical form is quite low due to absence of higher order MSIs \cite{kloeden1992stochastic,tripura2020ito}, since there is no suitable mechanism to treat the drift term appropriately. The inappropriate treatment results in the approximation error particularly due to inadequate attention on the non-linear part of the draft. This non-linear approximation error can be absorbed in the volatility component of the diffusion introducing a shift which can be treated using change of measure of Brownian motion. If the non-linear approximation error is $\mathcal{F}$-measurable then the change of measure due to this error can be formulated using Girsanov transformation. In this section it is shown how one can use the Girsanov change of measure to formulate a framework to treat the errors arising due to approximate treatment of the non-linear terms. Further, the use of rejection sampling to approximate the estimate of stochastic exponential integrals and to evaluate the Radon-Nikodym derivative arising from the change of measure is portrayed. With a little abuse of notations, let us assume that ${\bm{X}}$ and ${\bm{\dot X}}$ represents the displacement and velocity states of a second order multivariate dynamical system, which is different from the symbolic representation of $X$ in previous sections. Then, without loss of generality let us consider the $m$-dimensional oscillator purely excited by $n$-dimensional Brownian process:
\begin{equation}\label{3.1}
{\bf M}{\bm{\ddot X}} + {\bf C}{\bm{\dot X}} + {\bf{K}}{\bm X} + {\bf{\Omega }}\left( {t,{\bm{X}},{\bm{\dot X}}} \right) = \sum\limits_{k = 1}^n {{{\bf{F}}_k}(t,{\bm{X}},{\bm{\dot X}}){{\dot B}_k}(t)} 
\end{equation}
where, ${\bm X} = {\left\{ {{X_1},{X_2}, \cdots ,{X_m}} \right\}^T}$ is the response vector subjected to initial conditions ${\bm X}_{t_0} = \left\{ {{X_0},{X_0},\cdots,{X_0}} \right\}$, ${\bf{M}}\in {\mathbb{R} ^{m \times m}}$, ${\bf{C}}\in {\mathbb{R} ^{m \times m}}$ and ${\bf{K}}\in {\mathbb{R} ^{m \times m}}$ are respective constant mass, damping and stiffness matrices, ${\bf{\Omega }}\left( t,{{\bm{X}},{\bm{\dot X}}} \right)$ is a non-linear function of states (not necessarily smooth) but ensures a unique solution to Equation \ref{3.1} at least in the weak sense, $\left\{ {{{{\bf{F}}_k}}:\mathbb{R}  \mapsto {\mathbb{R} ^n}} \right\}$ is the set of n diffusion vectors (additive), ${B_k}$ for ${k \in [1,n]}$ is independently evolving zero-mean Brownian processes. The description of the oscillator as in Eq. (\ref{3.1}) is entirely formal (due to non-differentiability of the Brownian process, which implies that ${{{\dot B}_k}}$ exists merely as a valid measure, but not as a mathematical function). Introducing the state-space transformation ${\bm{X}} = {{\bm{X}}_1}$ and ${\bm{\dot X}} = {{\bm{X}}_2}$, the corresponding system of 2$m$-first order It\^{o}-SDEs for Eq. (\ref{3.1}) can be obtained as:
\begin{equation}\label{3.2}
\begin{array}{l}
d{X_{1j}}(t) = {X_{2j}}dt\\
d{X_{2j}}(t) = {g_j}\left( {t,{\bm X}} \right)dt + \sum\limits_{k = 1}^n {{f_{jk}}(t)d{B_k}} 
\end{array}
\end{equation}
where, ${X_{1j}}$ and ${X_{2j}}$ are the displacement and velocity states at the $j-(th)$ DOF for $j = \{1,2,3.......m\}$, $m$ being the number of DOF. Noting the normalization ${\bf{\tilde C}} = {{\bf{M}}^{ - 1}}{\bf{C}}$, ${\bf{\tilde K}} = {{\bf{M}}^{ - 1}}{\bf{K}}$ and ${\bf{\tilde \Omega }} = {{\bf{M}}^{ - 1}}{\bf{\Omega }}$, the drift term is identified as,
\begin{equation}\label{3.3}
{g_j}\left( {t,{\bm X}} \right) =  - \sum\limits_{k = 1}^m {{{\tilde C}_{jk}}{X_{2k}}}  - \sum\limits_{k = 1}^m {{{\tilde K}_{jk}}{X_{1k}} - \tilde{\Omega}_j \left( {t,{\bm X}} \right)} 
\end{equation}
The complete state vector is given by ${\bm{X}_j} = {\left\{ {X_{1j}},{X_{2j}} \right\}^T}$ and ${\sigma _{jk}}$ is the $(j,k)$-${th}$ element of diffusion matrix ${\bf{F}}$. Assuming that the drift vector ${g_j}\left( {t,{\bm X}} \right)$ can be decomposed into linear and non-linear parts as: ${g_j}\left( {t,{\bm X}} \right) = g_j^l\left( {t,{\bm X}} \right) + g_j^{nl}\left( {t,{\bm X}} \right)$, the Eq. (\ref{3.2}) can be rephrased as:
\begin{equation}\label{3.4}
\begin{array}{l}
d{X_{1j}}(t) = {X_{2j}}dt\\
d{X_{2j}}(t) = \left( {g_j^l\left( {t,{\bm X}} \right) + g_j^{nl}\left( {t,{\bm X}} \right)} \right)dt + \sum\limits_{k = 1}^n {{f_{jk}}(t)d{B_k}} 
\end{array}
\end{equation} 
The evolution of state ${X_{1j}}$ is comparatively straightforward as compared to that of ${X_{2j}}$, and does not require any special treatment. At this stage, one can approximate the non-linear function $g_j^{nl}\left( {t,{\bm X}} \right)$ over $t \in [{t_i},{t_{i - 1}}]$ by It\^{o}-Taylor expansion as:
\begin{equation}\label{3.5}
g_j^{nl}\left( {t,{\bm X}} \right) = g_j^{nl}\left( {{t_{i - 1}},{{\bm X}_{i - 1}}} \right) + \int_{{t_{i - 1}}}^{{t_i}} {{\Im ^0}\left( {g_j^{nl}\left( {s,{{\bm X}_s}} \right)} \right)} ds + \int_{{t_{i - 1}}}^{{t_i}} {{\Im ^1}\left( {g_j^{nl}\left( {s,{{\bm X}_s}} \right)} \right)} dB
\end{equation}
where, ${\Im ^0}(.)$ and ${\Im ^1}(.)$ are the SDE generators \cite{tripura2020ito}. Since, the aim this work is to keep the EM as simple and computationally efficient as its traditional form, by retaining the first term in the above expansion of $g_j^{nl}\left( {t,{\bm X}} \right)$, one can obtain the modified SDE of the Eq. (\ref{3.2}) as follows:
\begin{equation}\label{3.6}
\begin{array}{l}
d{X_{1j}}(t) = {X_{2j}}dt\\
d{X_{2j}}(t) = \left( {g_j^l\left( {t,{\bm X}} \right) + g_j^{nl}\left( {{t_{i - 1}},{{\bm X}_{i - 1}}} \right)} \right)dt + \sum\limits_{k = 1}^n {{f_{jk}}(t)d{B_k}} 
\end{array}
\end{equation}
A comparison between Eq. (\ref{3.4}) and Eq. (\ref{3.6}) provides an error process: ${\delta _j} = g_j^l\left( {t,{\bm X}} \right) - g_j^{nl}\left( {{t_{i - 1}},{{\bm X}_{i - 1}}} \right)$. Then, by incorporating the error process in Eq. (\ref{3.4}) a change of measure for Brownian process can be effected by modifying the $P$-Brownian motion ${B(t)}$ as follows,
\begin{equation}\label{3.7}
\begin{array}{ll}
d{X_{2j}}(t) &= \left( {g_j^l\left( {t,{\bm X}} \right) + g_j^{nl}\left( {{t_{i - 1}},{{\bm X}_{i - 1}}} \right) + {\delta _j}} \right)dt + \sum\limits_{k = 1}^n {{f_{jk}}(t)d{B_k}} \\
&= \left( {g_j^l\left( {t,{\bm X}} \right) + g_j^{nl}\left( {{t_{i - 1}},{{\bm X}_{i - 1}}} \right)} \right)dt + \sum\limits_{k = 1}^n {\left( {{\delta _j}dt + {f_{jk}}(t)d{B_k}} \right)} 
\end{array}
\end{equation}
Without any loss of generality, applying: $\sum_{k = 1}^n {\left( {{\delta _j}dt + {f_{jk}}(t)d{B_k}} \right)}  \to \sum_{k = 1}^n {{f_{jk}}(t)\left( {f_{jk}^{ - 1}(t){\delta _j}dt + d{B_k}} \right)} $ and substituting the appropriate error process, the following SDEs in the measure $P$ can be easily obtained.
\begin{equation}\label{3.12}
\begin{array}{l}
d{X_{1j}} = {X_{2j}}dt\\
d{X_{2j}} = \left( {g_j^l\left( {t,{\bm X}} \right) + g_j^{nl}\left( {{t_{i - 1}},{{\bm X}_{i - 1}}} \right)} \right)dt + \sum\limits_{k = 1}^n {{f_{jk}}(t)\left\{ {f_{jk}^{ - 1}(t)\left( {g_j^n\left( {t,{\bm X}} \right) - g_j^{nl}\left( {{t_{i - 1}},{{\bm X}_{i - 1}}} \right)} \right)dt + d{B_k}(t)} \right\}} 
\end{array}
\end{equation}
The main aim of this work is to simulate the system in an equivalent form where the error gets eliminated without effecting the dynamics of the system. This means that the present framework must remove the error so that the SDE in measure $Q$ takes the following form:  
\begin{equation}
d{X_{2j}} = \left( {g_j^l\left( {t,{\bm X}} \right) + g_j^{nl}\left( {{t_{i - 1}},{{\bm X}_{i - 1}}} \right)} \right)dt + \sum\limits_{k = 1}^n {{f_{jk}}(t)d{{\tilde B}_k}(t)} 
\end{equation}
Referring to Eq. (\ref{qbrownian}) the scalar valued shifting process, $\gamma \left( {t,{\bm X}} \right)$ can be identified from the integral representation of above equation as, $\gamma \left( {t,{\bm X}} \right) = \sum_{j = 1}^m {\sum_{k = 1}^n {f_{jk}^{ - 1}(t)\left( {g_j^l\left( {t,{\bm X}} \right) - g_j^{nl}\left( {{t_{i - 1}},{{\bm X}_{i - 1}}} \right)} \right)} } $. The $Q$-Brownian process ${{\tilde B}_k}(t)$ is formulated as:
\begin{equation}
d{{\tilde B}_k}(t) = \sum\limits_{k = 1}^n {f_{jk}^{ - 1}(t)\left( {g_j^l\left( {t,{\bm X}} \right) - g_j^{nl}\left( {{t_{i - 1}},{{\bm X}_{i - 1}}} \right)} \right)dt + d{B_k}(t)} 
\end{equation}
The change of measure $P \to Q$ is then constructed using the Radon-Nikodym derivative as follows:
\begin{equation}
\Lambda \left( {{t_i}} \right) = \exp \left( {\sum\limits_{k = 1}^n {\left( {\int_{{t_{i - 1}}}^{{t_i}} {{\gamma _k}\left( {s,{\bm X}} \right)d{{\tilde B}_k}(s)}  - \frac{1}{2}\int_{{t_{i - 1}}}^{{t_i}} {{{\left( {{\gamma _k}\left( {s,{\bm X}} \right)} \right)}^2}ds} } \right)} } \right)
\end{equation}
For ease of understanding, let the Radon-Nikodym derivative are computed as: $\Lambda ({t_i}) = \prod_{k = 1}^n {{\Lambda _k}({t_i})} $. Noting that ${\gamma _k}\left( {t,{\bm X}} \right) = \sum_{j = 1}^m {f_{jk}^{ - 1}(t){\delta _j}\left( {t,{{\bm X}_t},{{\bm X}_{i - 1}}} \right)} $ for $k=\{1, \ldots, n\} $, the term $\Lambda_k (t_i)$ is given as:
\begin{equation}\label{3.13}
\Lambda_k (t_i) = \exp \left( {\int_{{t_{i - 1}}}^{{t_i}} {{\gamma _k}\left( {s,{\bm X}} \right)} d{{\tilde B}_k} - \frac{1}{2}\int_{{t_{i - 1}}}^{{t_i}} {{{\left( {{\gamma _k}\left( {s,{\bm X}} \right)} \right)}^2}ds} } \right)
\end{equation}
The first integral term can be further simplified by expanding the stochastic exponential into a series of It\^{o} integrals. Towards this, let us exploit the It\^{o}-product rule (stochastic integration by parts):
\begin{equation}
d\left( {{\gamma _k}{{\tilde B}_k}} \right) = {\gamma _k}d{{\tilde B}_k} + {{\tilde B}_k}d{\gamma _k} + d\left[ {{\gamma _k}{{\tilde B}_k}} \right]
\end{equation}
Rearranging the terms, the integral representation of the above expression then follows:
\begin{equation}\label{3.13a}
\int_{{t_{i - 1}}}^{{t_i}} {{\gamma _k}} d{{\tilde B}_k} = \int_{{t_{i - 1}}}^{{t_i}} {d\left( {{\gamma _k}{{\tilde B}_k}} \right)}  - \int_{{t_{i - 1}}}^{{t_i}} {{{\tilde B}_k}d{\gamma _k}}  - \int_{{t_{i - 1}}}^{{t_i}} {d{\gamma _k}d{{\tilde B}_k}} 
\end{equation}
Here, the It\^{o} lemma for 2-dimensions with quadratic identities \cite{hassler2016stochastic}: $d{t^2} = 0$, $dB_t^2 = dt$ and $d{B_t}dt = 0$, can be applied to evaluate $d{\gamma _k}$ by incorporating the diffusions from Eq. (\ref{3.2}) as:
\begin{equation}\label{3.14}
d{\gamma _k} = {\left( {\frac{{\partial {\gamma _k}}}{{\partial {x_1}}}} \right)^T}d{x_1} + {\left( {\frac{{\partial {\gamma _k}}}{{\partial {x_2}}}} \right)^T}d{x_2} + \frac{1}{2}{\left( {\frac{{{\partial ^2}{\gamma _k}}}{{\partial x_2^2}}d{x_2}} \right)^T}d{x_2}
\end{equation}
Upon substitution of Eq. (\ref{3.14}) in to Eq. (\ref{3.13a}), the following is obtained (more details are available in \cite{raveendran2013nearly}):
\begin{equation}\label{3.15}
\begin{array}{ll}
\int\limits_{{t_{i - 1}}}^{{t_i}} {{\gamma _k}} d{{\tilde B}_k} =& \left[ {{\gamma _k}{{\tilde B}_k}} \right]_{{t_{i - 1}}}^{{t_i}} - \left[ {\tilde B_k^2{{\left( {{\partial _{{x_2}}}{\gamma _k}} \right)}^T}{f_k}} \right]_{{t_{i - 1}}}^{{t_i}} - \int\limits_{{t_{i - 1}}}^{{t_i}} {{{\left( {{\partial _{{x_2}}}{\gamma _k}} \right)}^T}{f_k}ds}  - \int\limits_{{t_{i - 1}}}^{{t_i}} {{{\tilde B}_k}{{\left( {{\partial _{{x_1}}}{\gamma _k}} \right)}^T}{x_2}ds} \\
& - \int\limits_{{t_{i - 1}}}^{{t_i}} {{{\tilde B}_k}{{\left( {{\partial _{{x_2}}}{\gamma _k}} \right)}^T}{{\tilde g}_j}ds}  - \frac{1}{2}\int\limits_{{t_{i - 1}}}^{{t_i}} {{{\tilde B}_k}{{\left( {\partial _{{x_2}}^2{\gamma _k}} \right)}^T}{f_k}ds}  + \int\limits_{{t_{i - 1}}}^{{t_i}} {{{\left( {\tilde B_k^2\left[ {\partial _{{x_1}{x_2}}^2{\gamma _k}} \right]{f_k}} \right)}^T}{x_2}ds} \\
& + \int\limits_{{t_{i - 1}}}^{{t_i}} {{{\left( {\tilde B_k^2\left[ {\partial _{{x_2}}^2{\gamma _k}} \right]{f_k}} \right)}^T}{{\tilde g}_j}ds}  + \int\limits_{{t_{i - 1}}}^{{t_i}} {{{\left( {{{\tilde B}_k}\left[ {\partial _{{x_2}}^2{\gamma _k}} \right]{f_k}} \right)}^T}{f_k}ds}  + \int\limits_{{t_{i - 1}}}^{{t_i}} {{{\left( {\tilde B_k^2\left[ {\partial _{{x_2}}^2{\gamma _k}} \right]{f_k}} \right)}^T}{f_k}d{{\tilde B}_k}} 
\end{array}
\end{equation}
Here, ${\partial _{{x_i}}}(.)$ and $\partial _{{x_i}}^2(.) $ are the first and second order partial derivatives with respect to state ${x_i}$. For the non-linear systems in structural and mechanical science, the second order derivative with respect to the state variable $x_2$ gets nullified, hence the It\^{o}-integral in  Radon-Nikodym derivative gets eliminated. For further treatment the boundary and integral terms can be separated and presented as ${\Lambda _k}({t_i}) = \Lambda _k^{(1)}({t_i})\Lambda _k^{(2)}({t_i})$, thus the Radon-Nikodym derivative can be now obtained as:
\begin{equation}\label{4.6}
\begin{array}{ll}
\Lambda ({t_i}) &= \prod\limits_{k = 1}^n {\Lambda _k^{(1)}({t_i})\Lambda _k^{(2)}({t_i})} \\
&= \prod\limits_{k = 1}^n {\underbrace {\exp \left( {\left[ {{\gamma _k}{{\tilde B}_k}} \right]_{{t_{i - 1}}}^{{t_i}} - \left[ {\tilde B_k^2{{\left( {{\partial _{{x_2}}}{\gamma _k}} \right)}^T}{f_k}} \right]_{{t_{i - 1}}}^{{t_i}}} \right)}_{\Lambda _k^{(1)}({t_i})}.\underbrace {\exp \left( { - \int_{{t_{i - 1}}}^{{t_i}} {{\phi _k}({{\bm X}_t},{{\bm X}_{i - 1}},{{\tilde B}_s})ds} } \right)}_{\Lambda _k^{(2)}({t_i})}} 
\end{array}
\end{equation}
Recalling, ${\gamma _k}\left( {t,{\bm X}} \right) = \sum_{j = 1}^m {f_{jk}^{ - 1}(t){\delta _j}\left( {t,{{\bm X}_t},{{\bm X}_{i - 1}}} \right)} $, the term ${{\phi _k}({{\bm X}_t},{{\bm X}_{i - 1}},{{\tilde B}_t})}$ is identified from Eq. (\ref{3.15}) as:
\begin{multline}
{\phi _k}({{\bm X}_t},{{\bm X}_{i - 1}},{{\tilde B}_t}) =  - {\left( {{\partial _{{x_2}}}{\gamma _k}} \right)^T}{f_k} - {{\tilde B}_k}{\left( {{\partial _{{x_1}}}{\gamma _k}} \right)^T}{X_2} - {{\tilde B}_k}{\left( {{\partial _{{x_2}}}{\gamma _k}} \right)^T}{{\tilde g}_j} - \frac{1}{2}{{\tilde B}_k}{\left( {\partial _{{x_2}}^2{\gamma _k}} \right)^T}{f_k} \\ + {\left( {\tilde B_k^2\left[ {\partial _{{x_1}{x_2}}^2{\gamma _k}} \right]{f_k}} \right)^T}{X_2} + {\left( {\tilde B_k^2\left[ {\partial _{{x_2}}^2{\gamma _k}} \right]{f_k}} \right)^T}{{\tilde g}_j} + {\left( {{{\tilde B}_k}\left[ {\partial _{{x_2}}^2{\gamma _k}} \right]{f_k}} \right)^T}{f_k} - \frac{1}{2}{\left( {{\gamma _k}\left( {s,{\bm X}} \right)} \right)^2}
\end{multline}
Although $\Lambda _k^{(2)}({t_i})$ does not explicitly contains stochastic terms, the integrals are still a function of Brownian process. The use of rejection sampling for evaluation of stochastic integrals reduce the complexity in such situation. The concept of rejection sampling for evaluation of $H({X_t}) = \int_{{t_{i - 1}}}^{{t_i}} {\phi (s,{B_s})ds} $ type of integrals can be pursued in detail in \cite{beskos2005exact}. Here, the idea is to form a binary indicator $\mathcal{I}$ such that $P \left( {\mathcal{I} = 1\left| {{\bm X}_t} \right.} \right) = \exp \left( { - \int_{{t_{i - 1}}}^{t_i} {\phi ds} } \right)$, where ${X_t}$ is finite collection of sample trajectory simulated between ${t_{i - 1}} \to {t_i}$. The particle is accepted with probability: $\exp \left( { - \int_{{t_{i - 1}}}^{t_i} {\phi ds} } \right)$, otherwise, rejected with probability: $1- \exp \left( { - \int_{{t_{i - 1}}}^{t_i} {\phi ds} } \right)$. To ensure that the error is bounded i.e. ${\left\| {{\delta _j}} \right\|_2} < \infty $, a constraint on ${\phi (t,{B_t})}$ is formed as:
\begin{equation}
0 \le \phi (t,{B_t}) \le \Delta _t^{ - 1}
\end{equation}
where, ${\Delta _t} = {t_i} - {t_{i - 1}}$ is the time step. In order to ensure that acceptance of particles are more, the particles having probability: $\exp \left( { - \int_{{t_{i - 1}}}^{t_i} {\phi ds} } \right)$ can be obtained by drawing random variables $\left( {u,v} \right) \sim \mathcal{U} \left( {\left[ {{t_{i - 1}},{t_i}} \right] \times \left( {0,\Delta _t^{ - 1}} \right)} \right)$ such that $\left\{ {\phi (u) \ge v} \right\}$ will have the desired probability. The RN-derivative $\Lambda (t)$ is a positive quantity ($>0$) and it can be easily shown that for $B(t_0)=0$ the RN-derivative $\Lambda (t_0)=1$. The aim of the rejection sampling is to sample a path such that the error ${\left\| {{\delta _j}} \right\|_2}$ is minimized thereby the process ${\left\| \gamma  \right\|_2}$ becomes infinitesimal small quantity. This automatically renders the RN-derivative $\Lambda (t)$ close to 1, ensuring a higher sample acceptance ratio in the rejection sampling. The indicator function $\mathcal{I}$ is constructed as $w \sim \mathcal{U}(0, 1)$. Then the occurrence of the event having probability: $1- \exp \left( { - \int_{{t_{i - 1}}}^{t_i} {\phi ds} } \right)$ suffices following condition \cite{beskos2005exact}: 
\begin{equation}
\phi \left( {X\left( u \right)} \right) < v \quad \text{or} \quad w > \frac{1}{{k!}}
\end{equation}
For the dicretization $({t_{i - 1}} = {t_1}) \le {t_2} \ldots  \le ({t_n} = {t_i})$, one can find an ordered sequence $\phi \left( {{\bm X} \left( {{t_r}} \right)} \right);{r \in [1,n]}$ for simulated uniform random numbers within $[{t_i},{t_{i - 1}}]$. Then with the aid of rejection sampling the accepted paths are determined. The accepted paths are then resampled using the weight $\Lambda _k^{(1)}({t_r})$ for the set $\left\{ {{t_{i - 1}} < {t_t} < {t_i}:r \in [1,n]} \right\}$, which completes the correction of the EM approximated samples paths. Finally, required resampling of the set of paths ${X({t_r})}$ are done in order to obtain the target moments. The EM sample paths for the SDE in Eq. (\ref{3.2}) over interval $t \in [{t_i},{t_{i - 1}}]$ are created using the following discretization:
\begin{equation}
\begin{array}{*{20}{l}}
{{X_{1j}}({t_i}) = {X_{1j}}({t_{i - 1}}) + {X_{2j}}({t_{i - 1}})dt}\\
{{X_{2j}}({t_i}) = {X_{2j}}({t_{i - 1}}) + \left( {g_j^n\left( {{t_{i - 1}},{{\bm X}_{i - 1}}} \right) + g_j^{nl}\left( {{t_{i - 1}},{{\bm X}_{i - 1}}} \right)} \right)dt + \sum\limits_{k = 1}^n {{f_{jk}}(t)d{B_k}} }
\end{array}
\end{equation} 
A pseudo code for the proposed framework, following the line of development in \cite{beskos2005exact} is provided in \ref{algvecvec}.
\begin{algorithm}
	\caption{Pseudo code for proposed near Exact Euler-Maruyama framework}\label{algvecvec}
	\begin{algorithmic}[1]
		\State $X(t_0)=0$, and form N partition for interval [0,T]
		\For {${t_i} \in [{t_1},{t_N}]$}	
		\State Generate a sequence within $t_i$:= ${t_r} \in \left\{ {{t_{i - 1}} = {t_1} \le {t_2} \ldots  \le {t_n} = {t_i}} \right\}$
		\State Draw w $\sim$ Unif(0, 1) for the binary indicator $\mathcal{I}$ 
		\State Set a sequence counter k = 0
		\For {$r \in [1,n]$, within the time instant $t_i$}
		\State Draw $\left( {u,v} \right) \sim \text{Unif}\left( {0,\Delta _{{t_i}}^{ - 1}} \right)$ \label{keystep1}
		\State Update the counter k = k+1
		\State Find $\Delta {B_r} = B({t_r}) - B({t_{r - 1}})$
		\State Find the Euler-Maruyama trajectory ${\bm X}(u({t_t}))$
		\State Evaluate $\phi \left( {X\left( {{u_{{t_r}}}} \right)} \right)$ at instant ${t_r}$
		\If {$\phi \left( {{\bm X}\left(  {{u_{{t_r}}}} \right)} \right) <  {{v_{{t_r}}}}$ or ${w} > \frac{1}{{k!}}$}
		\If {k is even,}
		\State $\phi \left( {X\left( {{u_{{t_r}}}} \right)} \right)$ is retained
		\State $\Lambda _k^{(2)}({t_i}) = \exp \left( { - \int_{{t_{i - 1}}}^{{t_i}} {{\phi _k}ds} } \right)$
		\Else
		\State {$\mathcal{I}$=0, and $\phi \left( {X\left( {{u_{{t_r}}}} \right)} \right)$ is rejected} \label{keystep2}
		\EndIf
		\Else
		\State Repeat Step. \ref{keystep1}$\to$\ref{keystep2}
		\EndIf
		\EndFor
		\State Perform the resampling using the weight $\Lambda _k^{(1)}({t_i})$ to the accepted particle sets $\left\{ {X({t_t}):r \in [1,n]} \right\}$
		\EndFor
	\end{algorithmic}
\end{algorithm}

\section{Numerical demonstrations:}
This section extends the understanding of the aforementioned framework and demonstrates application of near Exact Euler-Maruyama (n-EEM) method over few class of nonlinear oscillators. Three nonlinear systems have been taken: (i) A Rayleigh Duffing-Van der pol (RDVP) oscillator driven by additive stochastic noise, (ii) A Duffing-Van der pol (DVP) oscillator driven by multiplicative noise and deterministic sinusoidal excitation (iii) a 2-DOF non-linear system. To evaluate the improvement over the existing GCLM method, the solutions are compared with the available exact method like the FPK equation \cite{mamis2015exact} and more accurate approximate methods like It\^{o} Taylor weak 3.0 scheme \cite{tripura2020ito}.

\subsection{Rayleigh Van der pol (RVP) oscillator driven by additive stochastic noise}
The Rayleigh Van der pol oscillator is taken from the study \cite{mamis2015exact}, whose governing equation of motion is:
\begin{equation}\label{4.13}
{\ddot X}(t) + \left( {{h_1} + {h_3}{X^2}(t) + {h_3}{{\dot X}^2}(t)} \right){\dot X}(t) + X(t) = \sigma \dot B\left( t \right)
\end{equation}
where, ${h_1}$ and ${h_3}$ are the scalar valued parameters of the RVP oscillator, ${\sigma}$ is the noise intensity of the stochastic force ${ B\left( t \right)}$. The stochastic force is modeled as zero mean Gaussian white noise ($\mathcal{N} \sim (0,1)$). Using the state-space transformation ${X = {X_1}}$ and ${\dot X = {X_2}}$ Eq. (\ref{4.13}) can be rewritten in the form of It\^{o}-diffusion SDE as:
\begin{equation}\label{4.14}
\begin{array}{l}
d{X_1}(t) = {X_2}(t)dt\\
d{X_2}(t) = \left( {\left( { - {h_1} - {h_3}X_1^2(t) - {h_3}X_2^2(t)} \right){X_2}(t) - {X_1}(t)} \right)dt + \sigma dB\left( t \right)
\end{array}
\end{equation}
with respective drift and diffusion matrix as:
\begin{equation}
{\bf{g}}\left( {t,{{\bm X}_t}} \right) = \left[ {\begin{array}{*{20}{c}}
	{{X_2}(t)}\\
	{\left( { - {h_1} - {h_3}X_1^2(t) - {h_3}X_2^2(t)} \right){X_2}(t) - {X_1}(t)}
	\end{array}} \right],\qquad {\bf{f}}\left( {t,{{\bm X}_t}} \right) = \left[ {\begin{array}{*{20}{c}}
	0\\
	\sigma 
	\end{array}} \right]
\end{equation}
Let the time of integration be $t \in [0,T]$ and the interval is partitioned into N steps $0 = {t_0} <  \ldots  < {t_i} <  \ldots  < {t_N} = T$ for $i = 1 \ldots N$. If the increments are defined as $\Delta {B_i} = ({B_{i + 1}} - {B_i})$ and $\Delta {t_i} = ({t_{i + 1}} - {t_i})$, then one can generate the sample paths from the Euler Maruyama mapping of the above It\^{o}-SDEs as:
\begin{equation}
\begin{array}{l}
{X_1}\left( {{t_i}} \right) = {X_1}\left( {{t_{i - 1}}} \right) + {g_1}\left( {{t_{i - 1}},{{\bm X}_{{t_{i - 1}}}}} \right)\Delta t\\
{X_2}\left( {{t_i}} \right) = {X_2}\left( {{t_{i - 1}}} \right) + {g_2}\left( {{t_{i - 1}},{{\bm X}_{{t_{i - 1}}}}} \right)\Delta t + {f_2}\Delta {B_i}
\end{array}
\end{equation}
Considering the same form of It\^{o} diffusion equation as in Eq \ref{4.14}, the linear ${{\bf{g}}^l}$ and nonlinear ${{\bf{g}}^{nl}}$ coefficient of the drift term are obtained as:
\begin{equation}
{{\bf{g}}^l} = \left[ {\begin{array}{*{20}{c}}
	0&1\\
	{ - 1}&{ - {h_1}}
	\end{array}} \right], \qquad {{\bf{g}}^{nl}} = \left[ {\begin{array}{*{20}{c}}
	0\\
	{ - {h_3}X_1^2 - {h_3}{\dot X}_2^2}
	\end{array}} \right]
\end{equation}
Applying the Girsanov's transformation over the diffusion Eq \ref{4.14}, new diffusion equation with a different measure is obtained:
\begin{equation}\label{4.18}
\begin{array}{l}
d{X_1}(t) = {X_2}(t)dt\\
d{X_2}(t) = \left( { - {X_1}(t) - {h_1}{X_2}(t) - {h_3}X_{1,i - 1}^2{X_{2,i - 1}} - {h_3}X_{1,i - 1}^3} \right)dt + \sigma d\tilde B\left( t \right)
\end{array}
\end{equation}
with, $\delta = {h_3}X_1^2(t){X_2}(t) + {h_3}X_2^3(t) - {h_3}X_{1,i - 1}^2{X_{2,i - 1}} - {h_3}X_{2,i - 1}^3$. Here, $\tilde B\left( t \right)$ is the equivalent Brownian motion in $Q$-measure, defined as: 
\begin{equation}
d\tilde B\left( t \right) =  - {\sigma ^{ - 1}}\left( {{h_3}X_1^2(t){X_2}(t) + {h_3}X_2^3(t) - {h_3}X_{1,i - 1}^2{X_{2,i - 1}} - {h_3}X_{2,i - 1}^3} \right)dt + dB\left( t \right)
\end{equation}
Identifying, $\tilde g =  - {X_1}(t) - {h_1}{X_2}(t) - {h_3}X_{1,i - 1}^2{X_{2,i - 1}} - {h_3}X_{2,i - 1}^3$ and \\ $\gamma \left( {{{\bm X}_t},{{\bm X}_{i - 1}}} \right) =  - {\sigma ^{ - 1}}\left( {{h_3}X_1^2(t){X_2}(t) + {h_3}X_2^3(t) - {h_3}X_{1,i - 1}^2{X_{2,i - 1}} - {h_3}X_{2,i - 1}^3} \right)$, the associated Radon Nikodym derivative $\Lambda_t$ originating due to this change of measure $P \to Q$ can be obtained by substituting k=1 in Eq. (\ref{3.13}). As it is already discussed that the form of ${\Lambda _t}$ which being a probability density function in Eq. (\ref{3.13}) can not be obtained explicitly as it contains the stochastic exponential terms, needs to expanded using It\^{o}-Taylor expansion. After evaluation of the terms in the Eq. (\ref{3.15}), provided in \ref{filter1}, ${\Lambda _t}$ can be obtained as an exponential product of terms involving boundary and integral expressions as:
\begin{equation}
{\Lambda _{{t_i}}} = \underbrace {\exp \left( {\left. {\gamma \tilde B} \right|_{{t_{i - 1}}}^{{t_i}} - \left. {{{\tilde B}^2}\frac{{\partial \gamma }}{{\partial {X_2}}}\sigma } \right|_{{t_{i - 1}}}^{{t_i}}} \right)}_{{\Lambda _t}^{(1)}}.\underbrace {\exp \left( { - \int_{{t_{i - 1}}}^{{t_i}} {\phi \left( {{{\bm X}_s},{{\bm X}_{i - 1}},{{\tilde B}_s}} \right)ds} } \right)}_{{\Lambda _t}^{(2)}}
\end{equation}
where, ${\Lambda _t^{(1)}}$ and ${\Lambda _t^{(2)}}$ can be found as:
\begin{equation}
\begin{array}{l}
{\Lambda _{{t_i}}^{(1)}} = \exp \left[ \begin{array}{l}
- {\sigma ^{ - 1}}\left\{ \begin{array}{l}
\left( {{h_3}X_{1,i}^2{X_{2,i}} + {h_3}X_{2,i}^3} \right){{\tilde B}_i} - \left( {{h_3}X_{1,r}^2{X_{2,r}} + {h_3}X_{2,r}^3} \right){{\tilde B}_r}\\
- \left( {{h_3}X_{1,r}^2{X_{2,r}} + {h_3}X_{2,r}^3} \right)\left( {{{\tilde B}_i} - {{\tilde B}_r}} \right)
\end{array} \right\}\\
- {{\tilde B}_i}^2\left( {{h_3}X_{1,i}^2 + 3{h_3}X_{2,i}^2} \right) + {{\tilde B}_r}^2\left( {{h_3}X_{1,r}^2 + 3{h_3}X_{2,r}^2} \right)
\end{array} \right]\\
\Lambda _{{t_i}}^{(2)} = \exp \left\{ { - \int {\phi \left( s \right)dt} } \right\}
\end{array}
\end{equation}
Here, ${B_r} = {B_{{t_r}}}$, ${X_{1,i}} = {X_1}({t_i})$, ${X_{1,r}} = {X_1}({t_r})$ and ${B_i} \approx \left( {{B_i} - {B_r}} \right) \sim N\left( {0,\sqrt {{t_i} - {t_r}} } \right)$ are the re-sampled Brownian increments, ${t_r} = ({t_{i - 1}} \ldots {t_i})$ is the re-sampled time instants. The term $\phi \left( t \right)$ at each time instant $t_i$ is evaluated as:
\begin{multline}
\psi \left( t \right) =  - \left( {{h_3}X_{1,i}^2 + 3{h_3}X_{2,i}^2} \right) - {{\tilde B}_t}{\sigma ^{ - 1}}\left( {2{h_3}{X_{1,i}}{X_{2,i}}} \right){X_{2,i}} - {{\tilde B}_t}{\sigma ^{ - 1}}\left( {{h_3}X_{1,i}^2 + 3{h_3}X_{2,i}^2} \right)\tilde g \\ - \frac{1}{2}{{\tilde B}_t}\left( {6{h_3}{X_{2,i}}} \right)\sigma  + {{\tilde B}_r}^2\left( {2{h_3}{X_{1,i}}} \right){X_{2,i}} + {{\tilde B}_t}^2\left( {6{h_3}{X_{2,i}}} \right)\tilde gdt + {{\tilde B}_t}\left( { - 6{h_3}{X_{2,i}}} \right)\sigma  + \frac{1}{2}{\gamma ^2}
\end{multline}
The exact stationary density function for Rayleigh Van der pol oscillator can be found in \cite{mamis2015exact} as:
\begin{equation}
p\left( {{X_1},{X_2}} \right) = C\exp \left\{ { - \frac{{{\eta _3}}}{2}{X_2}^4 - \left( {{\eta _1} + 2{\eta _3}U\left( {{X_1}} \right)} \right)\left( {{\eta _2}^2 + U\left( {{X_1}} \right)} \right) - {\eta _1}U\left( {{X_1}} \right)} \right\}
\end{equation}
with following terms,
\[{{\eta _3} = \frac{{{h_3}}}{{2{\sigma ^2}{D_{11}}}} > 0}, \quad {{\eta _1} = \frac{{{h_1} + 2{h_3}U\left( {{X_1}} \right)}}{{2{\sigma ^2}{D_{11}}}}}, \quad {U\left( {{X_1}} \right) = \frac{{{X_1}^2}}{2}}\]
The value of constant ${C}$ is obtained by the property that ${\int {\int_{ - \infty }^\infty  {p\left( {{X_1},{X_2}} \right)d{X_1}d{X_2}} }  = 1}$. Once the density is available, the mean squared moments of system responses can be found with no effort.

\begin{figure}[H]
	\centering
	\includegraphics[width=\textwidth]{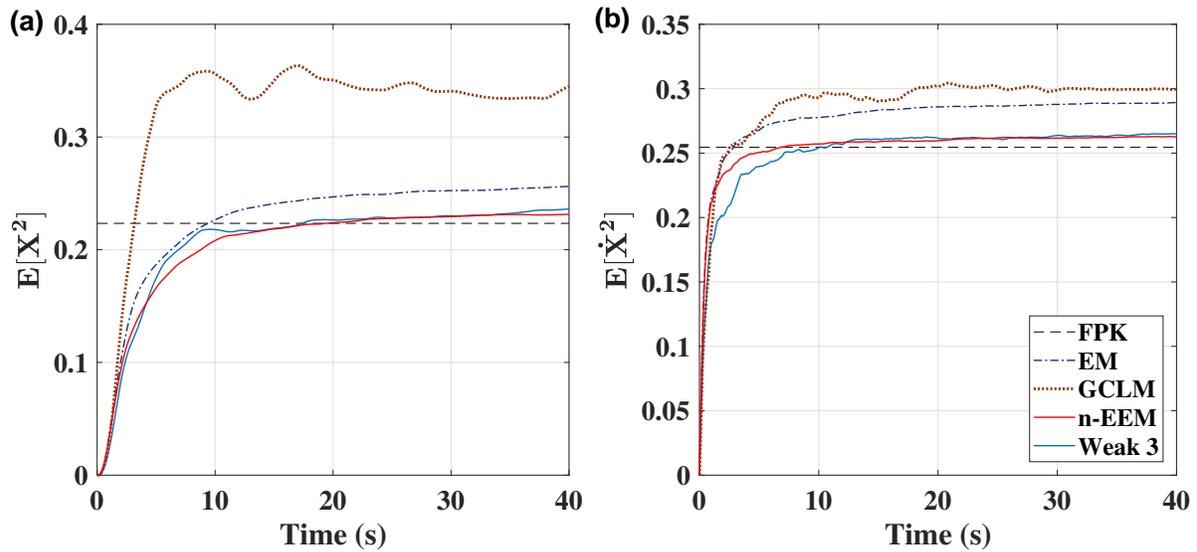}
	\caption{Sample averaged mean squared moments of displacement and velocity state of the RVP oscillator for ${h_1=1}$, ${h_3 = 1}$ and ${\sigma = 1}$}
	\label{fig:rdvpd}
\end{figure}

\begin{figure}[H]
	\centering
	\includegraphics[width=0.8\textwidth]{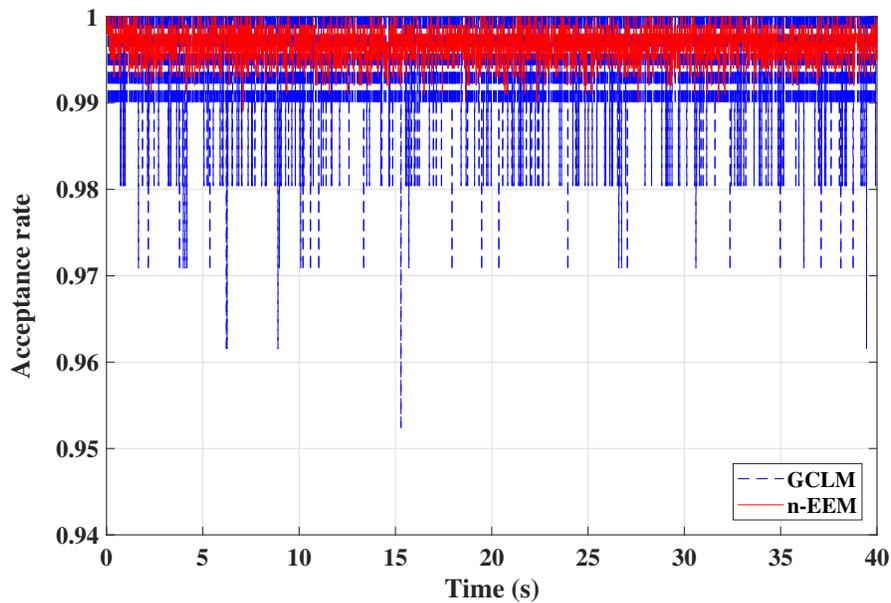}
	\caption{Acceptance ratio of the sample paths in rejection sampling for RVP oscillator for ${h_1=1}$, ${h_3 = 1}$ and ${\sigma = 1}$}
	\label{fig:rdvpv}
\end{figure}

The mean time responses of $E\left[ {{X^2}} \right]$ and $E\left[ {{{\dot X}^2}} \right]$ are computed using ${\Delta t = 0.01}$s with an ensemble size of 100. An initial condition of $\left( {{X_0},{{\dot X}_0}} \right) = \left( {0.01,0.01} \right)$ is assumed for the simulation. The parameters of the oscillator were taken as: ${h_1=1 }$, ${h_3 = 1}$ and ${\sigma = 1}$. The squared moments of system states along with the acceptance ratio from rejection sampling are portrayed in Fig. \ref{fig:rdvpd} and Fig. \ref{fig:rdvpv}. It is evident from Fig. \ref{fig:rdvpd}  that, as the time history progresses the second moment trajectory obtained using proposed n-EEM scheme approximates the stationary solution, whereas, the moment trajectories obtained from the EM and GCLM diverges far from the stationary solution. In Fig. \ref{fig:rdvpv}, the acceptance ratios obtained using GCLM and the proposed n-EEM are shown. It is noted that the proposed approach yields a higher acceptance ratio with comparatively less fluctuation; this is because the 
EM based proposal density used in this paper and the target density resides in close proximity of each other. On the contrary, the proposal density used in GCLM is a crude approximation of the target density and hence, fluctuations in acceptance ratio is observed. This comparison of the acceptance ratio for the proposed n-EEM and available GCLM scheme demonstrates the computational efficiency of the proposed scheme over available GCLM method that also shares a similar rejection sampling scheme.

\subsection{Duffing-Van der pol (DVP) oscillator driven by stochastic and sinusoidal excitation}
The system is taken from the literature \cite{tripura2020ito} where the non-linearity arises due to an additive cubic dissipation force. The dynamics of the system is governed by the following equation:
\begin{equation}\label{4.26}
m\ddot X\left( t \right) + c\dot X\left( t \right) - kX\left( t \right) + \alpha {X^3}\left( t \right) = \rho X\left( t \right)\dot B\left( t \right) + A\sin \left( {2\pi \omega t} \right)
\end{equation}
where, m, k and c are mass, stiffness and damping of the system, respectively, ${\alpha}$ is the parameter of the DVP oscillator, ${\rho}$ is the noise intensity of the stochastic force ${\dot B\left( t \right)}$. The force ${B\left( t \right)}$ is modeled as a zero mean Gaussian white noise and ${\omega}$ is the frequency of sinusoidal excitation. Through the state-space transformation ${X = {X_1}}$ and ${\dot X = {X_2}}$ Eq. (\ref{4.26}) can be rewritten in the form of It\^{o}-diffusion SDE as:
\begin{equation}\label{4.27}
d\left[ {\begin{array}{*{20}{c}}
	{{X_1}(t)}\\
	{{X_2}(t)}
	\end{array}} \right] = \underbrace {\left[ {\begin{array}{*{20}{c}}
		{{X_2}(t)}\\
		{ - \frac{1}{m}\left( { - k{X_1}(t) + c{X_2}(t) + \alpha {X_1}^3(t)} \right) + \frac{A}{m}\sin \left( {2\pi \omega t} \right)}
		\end{array}} \right]}_{{\bf{g}}\left( {t,{{\bm X}_t}} \right)}dt + \underbrace {\left[ {\begin{array}{*{20}{c}}
		0\\
		{\frac{{\rho {X_1}(t)}}{m}}
		\end{array}} \right]}_{{\bf{f}}\left( {t,{{\bm X}_t}} \right)}dB\left( t \right)
\end{equation}
One can find the Euler Maruyama mapping for the system follows:
\begin{equation}
\begin{array}{l}
{X_1}\left( {{t_i}} \right) = {X_1}\left( {{t_{i - 1}}} \right) + {g_1}\left( {{t_{i - 1}},{{\bm X}_{{t_{i - 1}}}}} \right)\Delta t\\
{X_2}\left( {{t_i}} \right) = {X_2}\left( {{t_{i - 1}}} \right) + {g_2}\left( {{t_{i - 1}},{{\bm X}_{{t_{i - 1}}}}} \right)\Delta t + {f_2}\left( {{t_{i - 1}},{{\bm X}_{{t_{i - 1}}}}} \right)\Delta {B_i}
\end{array}
\end{equation}
The linear and nonlinear drift coefficient matrices: ${{\bf{g}}^l}$ and ${{\bf{g}}^{nl}}$ are identified from the drift matrix ${\bf g} \left( {t,{{\bm X}_t}} \right)$ as:
\begin{equation}
{{\bf{g}}^l} = \left[ {\begin{array}{*{20}{c}}
	0&1\\
	{\frac{k}{m}}&{ - \frac{c}{m}}
	\end{array}} \right], \quad {{\bf{g}}^{nl}} = \left[ {\begin{array}{*{20}{c}}
	0\\
	{ - \frac{{\alpha {X_1}^3(t)}}{m}}
	\end{array}} \right]
\end{equation}
An equivalent SDE of the diffusion in Eq. (\ref{4.27}) in the $Q$-measure through Girsanov transformation is obtained as:
\begin{equation}
\begin{array}{l}
d{X_1}(t) = {X_2}(t)dt\\
d{X_2}(t) = \left( {\frac{k}{m}{X_1}(t) - \frac{c}{m}{X_2}(t) - \frac{\alpha }{m}X_{1,i - 1}^3} \right)dt + \frac{{\rho {X_1}(t)}}{m}\left[ { - {{\left( {\frac{{\rho {X_1}(t)}}{m}} \right)}^{ - 1}}\left( {\frac{\alpha }{m}X_1^3(t) - \frac{\alpha }{m}X_{1,i - 1}^3} \right)dt + \frac{{{X_1}(t)}}{m}dB\left( t \right)} \right]
\end{array}
\end{equation}
One can verify that the $Q$-Brownian motion is given as: 
\begin{equation}
d\tilde B\left( t \right) =  - {\left( {\frac{{\rho {X_1}(t)}}{m}} \right)^{ - 1}}\left( {\frac{\alpha }{m}X_1^3(t) - \frac{\alpha }{m}X_{1,i - 1}^3} \right)dt + dB\left( t \right)
\end{equation}
Further, noting that $\tilde g = \left(\frac{k}{m}{X_1}(t) - \frac{c}{m}{X_2}(t) - \frac{\alpha }{m}X_{1,i - 1}^3 \right)$ and $\gamma \left( {{{\bm X}_t},{{\bm X}_{i - 1}}} \right) =  - \left[{\left( {\frac{{\rho {X_1}(t)}}{m}} \right)^{ - 1}}\left( {\frac{\alpha }{m}X_1^3(t) - \frac{\alpha }{m}X_{1,i - 1}^3} \right) \right]$, the ${\Lambda _t}$ associated with this change of measure has the form in Eq. (\ref{3.13}). The detailed calculation for the estimate of the term: $\exp \left\{ {\int_{{t_{i - 1}}}^{{t_i}} {\gamma \left( {{{\bm X}_s},{{\bm X}_{i - 1}}} \right)d\tilde B\left( s \right)} } \right\}$ is given in \ref{filter2}. Noting that ${\Lambda _{{t_i}}}$ can be reduced in to the form: ${\Lambda _{{t_i}}} = \Lambda _{{t_i}}^{(1)}\Lambda _{{t_i}}^{(2)}$, they are found as follows:
\begin{equation}
\begin{array}{l}
\Lambda _{{t_i}}^{(1)} = \exp \left[ { - {{\left( {\frac{{\rho {X_1}(t)}}{m}} \right)}^{ - 1}}\left( {\frac{\alpha }{m}X_{1,i}^3{{\tilde B}_i} - \frac{\alpha }{m}X_{1,r}^3{{\tilde B}_r} - \frac{\alpha }{m}X_{1,i - 1}^3\left( {{{\tilde B}_i} - {{\tilde B}_r}} \right)} \right)} \right]\\
\Lambda _{{t_i}}^{(2)} = \exp \left\{ { - \int\limits_{{t_{i - 1}}}^{{t_i}} {\left( {{{\tilde B}_s}\left( { - 3{{\left( {\frac{{\rho {X_1}(s)}}{m}} \right)}^{ - 1}}\frac{\alpha }{m}X_1^2(s)} \right){X_2}(s) + \frac{1}{2}{\gamma ^2}} \right)dt} } \right\}
\end{array}
\end{equation} 

\begin{figure}[H]
	\centering
	\includegraphics[width=\textwidth]{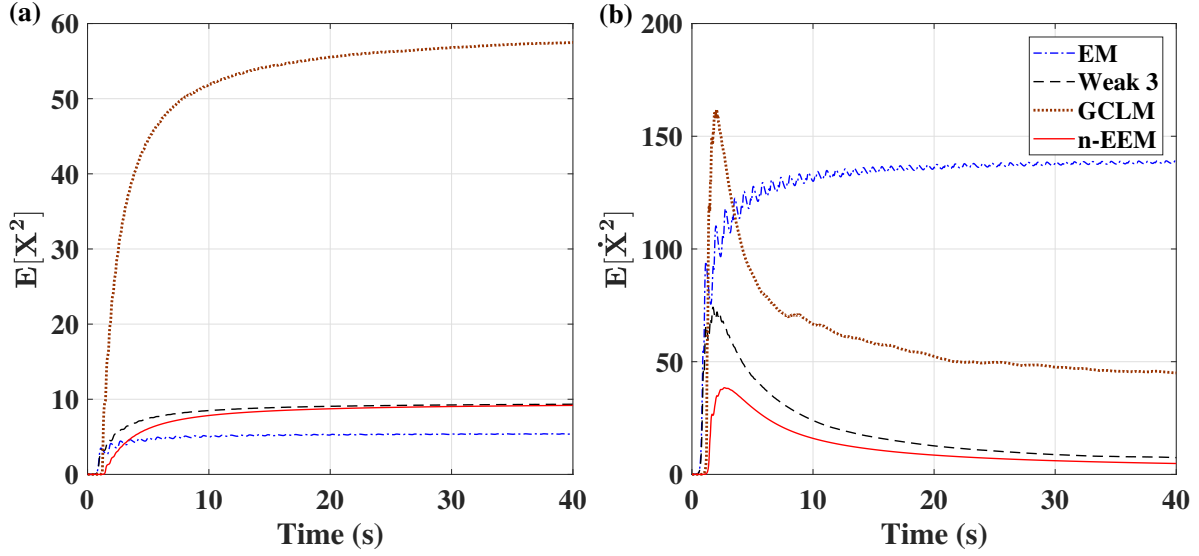}
	\caption{Mean squared responses of displacement and velocity states of DVP oscillator for ${m=1}$kg, ${k=60}$N/m, ${\alpha=2}$, ${c=1.5}$Ns/m and ${\rho=0.5}$}
	\label{fig:dvpd}
\end{figure}

\begin{figure}[H]
	\centering
	\includegraphics[width=0.8\textwidth]{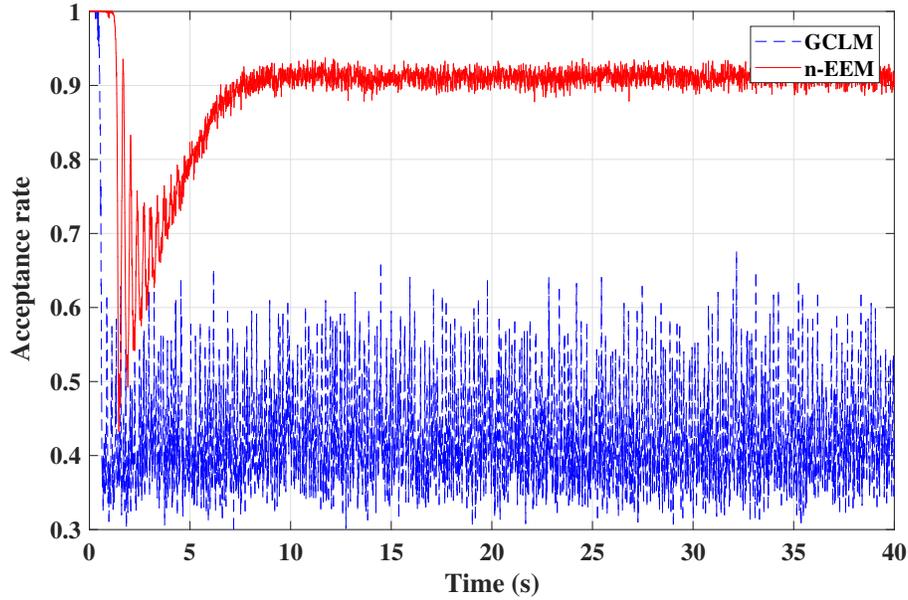}
	\caption{Sample acceptance ratio from rejection sampling for DVP oscillator ${m=1}$kg, ${k=60}$N/m, ${\alpha=2}$, ${c=1.5}$Ns/m and ${\rho=0.5}$}
	\label{fig:dvpv}
\end{figure}

The second moments $E\left[ {{X^2}} \right]$ and $E\left[ {{{\dot X}^2}} \right]$ for the system are computed using ${\Delta t = 0.01}$s with an effective ensemble size of 100. The parameters were taken as, ${m=}$1kg, ${k=60}$N/m, ${\alpha=2}$, ${c=}$1.5kNs/m and ${\rho=0.5}$. The initial conditions were: $\left( {{X_0},{{\dot X}_0}} \right) = \left( {0.01,0.01} \right)$. In absence of the exact stationary solution, the reference second-moment time histories of system displacement and velocity are generated using the weak 3.0 It\^{o}-Taylor scheme with 1000 Monte Carlo. The second order moment of the responses for this case study is presented in Fig. \ref{fig:dvpd}, which clearly shows that, as the time progresses the second moments obtained using the proposed n-EEM scheme converges to the reference solution, whereas, the second moments of the EM and GCLM solutions does not converges to the reference solution. In Fig. \ref{fig:dvpv} the acceptance ratios in proposed n-EEM is compared with the available GCLM method. initially, a relatively lower acceptance ratio is observed for n-EEM. This is because the EM based proposal density yields a poor approximation of the target density. However, as time progresses, the quality of the EM based proposal density improves as evident from the higher acceptance ratio. In case of GCLM method, the lower acceptance ratio throughout the time indicates its inefficiency.

\subsection{A 2-DOF non-linear system driven by stochastic excitation}
\begin{figure}[H]
	\centering
	\includegraphics[width=\textwidth]{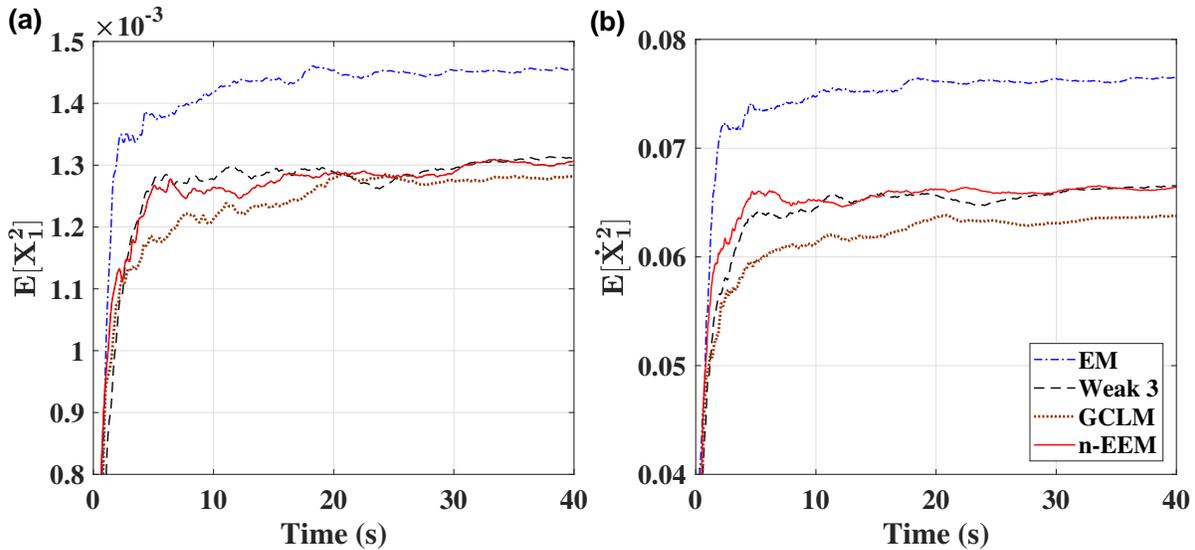}
	\caption{Mean squared histories of the displacement and velocity states of the 2-DOF oscillator (first DOF) with parameters: ${m_1 = m_2 =1}$kg, ${k_1 = k_2 =100}$N/m, ${c_1 = c_2 =7.75}$Ns/m, ${\alpha = \beta =100}$ and ${\sigma_1 = \sigma_2 = 1}$}
	\label{fig:2dofd11}
\end{figure}

\begin{figure}[H]
	\centering
	\includegraphics[width=\textwidth]{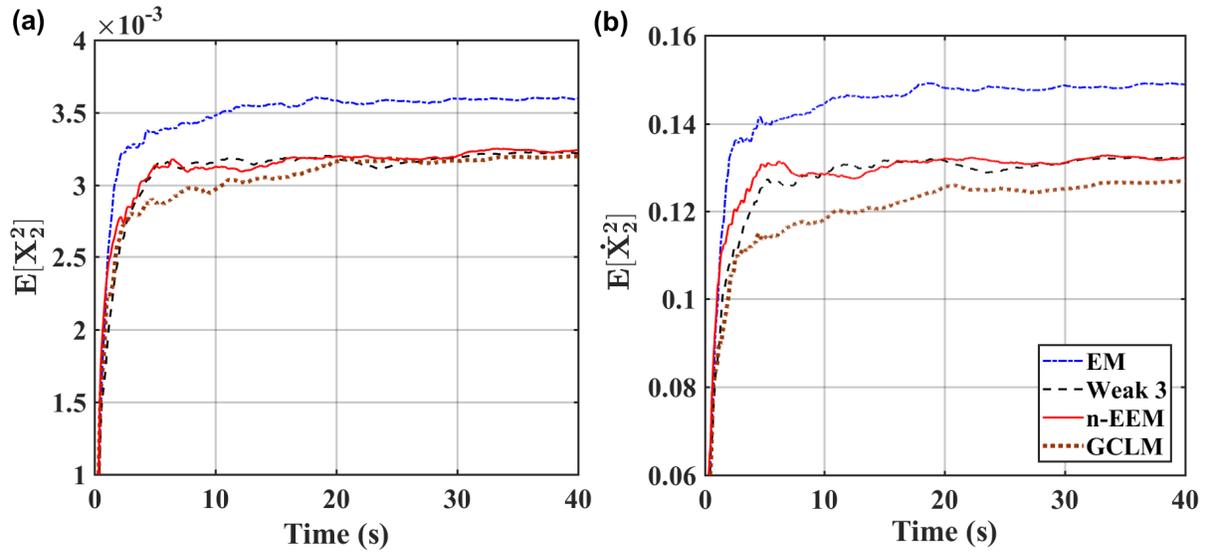}
	\caption{Mean squared histories of the displacement and velocity states of the 2-DOF oscillator (second DOF) with parameters: ${m_1 = m_2 =1}$kg, ${k_1 = k_2 =100}$N/m, ${c_1 = c_2 =7.75}$Ns/m, ${\alpha = \beta =100}$ and ${\sigma_1 = \sigma_2 = 1}$}
	\label{fig:2dofd22}
\end{figure}

This example demonstrates the possible extension of the proposed method towards an MDOF system. A 2-DOF non-linear oscillator is considered here whose dynamical equation of motions are sa follows:
\begin{equation}
\begin{array}{l}
{{\ddot X}_1}(t) + \left( {{c_1} + {c_2}} \right){{\dot X}_1}(t) - {c_2}{{\dot X}_2}(t) + \left( {{k_1} + {k_1}} \right){X_1}(t) - {k_2}{X_2}(t) + \alpha X_1^2(t){{\dot X}_1}(t) = {\sigma _1}{{\dot B}_1}\left( t \right)\\
{{\ddot X}_2}(t) - {c_2}{{\dot X}_1}(t) + {c_2}{{\dot X}_2}(t) - {k_2}{X_1}(t) + {k_2}{X_2}(t) + \beta X_2^3(t) = {\sigma _2}{{\dot B}_2}\left( t \right)
\end{array}
\end{equation} 
where, ${c_i}$ and ${k_i}$ are the respective damping and stiffness of the $i^{th}$ DOF of the oscillator, ${\alpha}$ and $\beta$ are the co-efficient  of non-linear dissipation forces, $\sigma_i$ is the noise intensity of the Brownian motion ${{{\dot B}_i}\left( t \right)}$ at $i^{th}$-DOF. A state-space model containing four variables: ${X_1} = {Y_1},{{\dot X}_1} = {Y_2},{X_2} = {Y_3},{{\dot X}_2} = {Y_4}$, helps to construct the first order It\^{o}-SDEs for the system in the following form: $d{\bm{Y}}(t) = {\bf{g}}\left( {t,{{\bm{Y}}_t}} \right)dt + {\bf{f}}\left( {t,{{\bm{Y}}_t}} \right)d{\bm{B}}(t)$, where, the respective drift and diffusion matrices for the system is:
\begin{equation}
\begin{array}{l}
{\bf{g}}\left( {t,{{\bm Y}_t}} \right) = \left[ {\begin{array}{*{20}{c}}
	{{Y_2}(t)}\\
	{ - \left\{ {\left( {{c_1} + {c_2}} \right){Y_2}(t) - {c_2}{Y_4}(t) + \left( {{k_1} + {k_1}} \right){Y_1}(t) - {k_2}{Y_3}(t) + \alpha Y_1^2(t){Y_2}(t)} \right\}}\\
	{{Y_4}(t)}\\
	{ - \left\{ { - {c_2}{Y_2}(t) + {c_2}{Y_4}(t) - {k_2}{Y_1}(t) + {k_2}{Y_3}(t) + \beta Y_3^3(t)} \right\}}
	\end{array}} \right] \\
{\bf{f}}\left( {t,{{\bm Y}_t}} \right) = \left[ {\begin{array}{*{20}{c}}
	0&0\\
	{{\sigma _1}}&0\\
	0&0\\
	0&{{\sigma _2}}
	\end{array}} \right]
\end{array}
\end{equation}
One can find the Euler Maruyama mapping for the above SDE as:
\begin{equation}
{\bm{Y}}({t_i}) = {\bm{Y}}({t_{i - 1}}) + {\bf{g}}\left( {{t_{i - 1}},{{\bm{Y}}_{{t_{i - 1}}}}} \right)\Delta {t_i} + {\bf{f}}\left( {{t_{i - 1}},{{\bm{Y}}_{{t_{i - 1}}}}} \right)\Delta {\bm{B}}({t_i})
\end{equation}
For the proposed method the drift matrix can be decomposed into linear and nonlinear matrices as:
\begin{equation}
{{\bf{g}}^l} = \left[ {\begin{array}{*{20}{c}}
	0&1&0&0\\
	{ - ({k_1} + {k_2})}&{ - ({c_1} + {c_2})}&{{k_2}}&{{c_2}}\\
	0&0&0&1\\
	{{k_2}}&{{c_2}}&{ - {k_2}}&{ - {c_2}}
	\end{array}} \right], \quad {{\bf{g}}^{nl}} = \left[ {\begin{array}{*{20}{c}}
	0\\
	{ - \alpha Y_1^2(t){Y_2}(t)}\\
	0\\
	{ - \beta Y_3^3(t)}
	\end{array}} \right]
\end{equation}
The transformed system of SDEs in $Q$-measure is then given by,
\begin{equation}
\begin{array}{l}
d{Y_1}(t) = {Y_2}(t)dt\\
d{Y_2}(t) =  - \left\{ {\left( {{c_1} + {c_2}} \right){Y_2}(t) - {c_2}{Y_4}(t) + \left( {{k_1} + {k_1}} \right){Y_1}(t) - {k_2}{Y_3}(t) + \alpha Y_1^2(t){Y_2}(t)} \right\}dt + {\sigma _1}d{{\tilde B}_1}\left( t \right)\\
d{Y_3}(t) = {Y_4}(t)dt\\
d{Y_4}(t) =  - \left\{ { - {c_2}{Y_2}(t) + {c_2}{Y_4}(t) - {k_2}{Y_1}(t) + {k_2}{Y_3}(t) + \beta Y_3^3(t)} \right\}dt + {\sigma _2}d{{\tilde B}_2}\left( t \right)
\end{array}
\end{equation}
where, the $Q$-Brownian motions are given as:
\begin{equation}
\begin{array}{l}
d{{\tilde B}_1}\left( t \right) =  - \sigma _1^{ - 1}\alpha \left( {Y_1^2(t){Y_2}(t) - Y_1^2({t_{i - 1}}){Y_2}({t_{i - 1}})} \right)dt + d{B_1}\left( t \right)\\
d{{\tilde B}_1}\left( t \right) =  - \sigma _2^{ - 1}\beta \left( {Y_3^3(t) - Y_3^3({t_{i - 1}})} \right)dt + d{B_2}\left( t \right)
\end{array}
\end{equation}
Identified that ${\gamma _1}\left( {{{\bm Y}_t},{{\bm Y}_{i - 1}}} \right) =  - \sigma _1^{ - 1}\alpha \left( {Y_1^2(t){Y_2}(t) - Y_1^2({t_{i - 1}}){Y_2}({t_{i - 1}})} \right)$ and ${\gamma _2}\left( {{{\bm Y}_t},{{\bm Y}_{i - 1}}} \right) =  - \sigma _2^{ - 1}\beta \left( {Y_3^3(t) - Y_3^3({t_{i - 1}})} \right)$, the Radon Nikodym derivative associated with this change of measure for the 2-DOF system can be verified as:
\begin{equation}\label{4.57}
{\Lambda _{{t_i}}} = \exp \left\{ {\int\limits_{{t_{i - 1}}}^{{t_i}} {{\gamma _1}\left( {{{\bm Y}_s},{{\bm Y}_{i - 1}}} \right)d{{\tilde B}_1}\left( s \right) + \int\limits_{{t_{i - 1}}}^{{t_i}} {{\gamma _2}\left( {{{\bm Y}_s},{{\bm Y}_{i - 1}}} \right)d{{\tilde B}_2}\left( s \right)} }  - \frac{1}{2}\int\limits_{{t_{i - 1}}}^{{t_i}} {\left( {{\gamma _1}{{\left( {{{\bm Y}_s},{{\bm Y}_{i - 1}}} \right)}^2} + {\gamma _2}{{\left( {{{\bm Y}_s},{{\bm Y}_{i - 1}}} \right)}^2}} \right)ds} } \right\}
\end{equation}
The evaluation of the stochastic exponential: $\exp \left\{ {\int_{{t_{i - 1}}}^{{t_i}} {{\gamma _2}\left( {{{\bm Y}_s},{{\bm Y}_{i - 1}}} \right)d{{\tilde B}_2}\left( s \right)} } \right\}$ and $\exp \left\{ {\int_{{t_{i - 1}}}^{{t_i}} {{\gamma _1}\left( {{{\bm Y}_s},{{\bm Y}_{i - 1}}} \right)d{{\tilde B}_1}\left( s \right)} } \right\}$ are provided in \ref{filter3}. The reduced form of $\Lambda (t)$ can then be obtained as:
\begin{equation}
{\Lambda _{{t_i}}} = \underbrace {\exp \sum\limits_{k = 1}^2 {\left( {\left. {{\gamma _k}{{\tilde B}_k}} \right|_{{t_{i - 1}}}^{{t_i}} - \left. {\tilde B_k^2\frac{{\partial {\gamma _k}}}{{\partial {y_{2k}}}}{\sigma _k}} \right|_{{t_{i - 1}}}^{{t_i}}} \right)} }_{\Lambda _{{t_i}}^{(1)}}\underbrace {\exp \left( { - \int_{{t_{i - 1}}}^{{t_i}} {\sum\limits_{k = 1}^2 {{\phi _k}({{\bm Y}_s},{{\bm Y}_{i - 1}},{{\tilde B}_s})} ds} } \right)}_{\Lambda _{{t_i}}^{(2)}}
\end{equation}
where, ${\varepsilon _1}(t) = Y_1^2(t){Y_2}(t) - Y_1^2({t_{i - 1}}){Y_2}({t_{i - 1}})$ and ${\varepsilon _2}(t) = Y_3^3(t) - Y_3^3({t_{i - 1}})$. Further, ${\Lambda _{{t_i}}^{(1)}}$ and ${\Lambda _{{t_i}}^{(2)}}$
\begin{equation}
\begin{array}{ll}
\Lambda _{{t_i}}^{(1)} =& \exp \Bigl\{  - \sigma _1^{ - 1}\alpha \left( {Y_{1,i}^2{Y_{2,i}}{{\tilde B}_{1,i}} - Y_{1,r}^2{Y_{2,r}}{{\tilde B}_{1,r}} - Y_{1,r}^2{Y_{2,r}}\left( {{{\tilde B}_{1,i}} - {{\tilde B}_{1,r}}} \right)} \right) \\&- \sigma _2^{ - 1}\beta \left( {Y_{3,i}^3 - Y_{3,r}^3 - Y_{3,r}^3\left( {{{\tilde B}_{1,i}} - {{\tilde B}_{1,r}}} \right)} \right) \Bigr\}\\
\Lambda _{{t_i}}^{(2)} =& \exp \left\{ { - \int {\psi \left( t \right)dt} } \right\}, \text{with}\\
\phi (t) =&  - \alpha Y_1^2(t) - 2{{\tilde B}_1}(t)\sigma _1^{ - 1}\alpha {Y_1}(t)Y_2^2(t) - {{\tilde B}_1}(t)\sigma _1^{ - 1}\alpha Y_1^2(t)\tilde g - 2B_1^2(t)\alpha {Y_1}(t){Y_2}(t) - 3\sigma _2^{ - 1}\beta Y_3^2(t){Y_4}(t)
\end{array}
\end{equation}
The simulation results for the 2-DOF oscillator is obtained using a Monte-Carlo ensemble of size 1000 at a time rate of $\Delta t=0.01$s. The system parameters are taken as: ${m_1 = m_2 =1}$kg, ${k_1 = k_2 =100}$N/m, ${c_1 = c_2 =7.75}$Ns/m, ${\alpha = \beta =100}$ and ${\sigma_1 = \sigma_2 = 1}$. An initial condition as: $\left( {{{\bf{X}}_0}} \right) = \left( {0.01,0.01,0.01,0.01} \right)$ is applied to the system. The reference solution obtained using higher order weak 3 It\^{o}-Taylor method. The second moment time histories of the first-DOF: $E\left[ {{X_1^2}} \right]$ and $E\left[ {{{\dot X}_1^2}} \right]$ are portrayed in Fig. \ref{fig:2dofd11} and that of the second-DOF: $E\left[ {{X_2^2}} \right]$ and $E\left[ {{{\dot X}_2^2}} \right]$ are presented in Fig. \ref{fig:2dofd22}. It can be observed from Figs. \ref{fig:2dofd11} and \ref{fig:2dofd22} that as the time progresses, the proposed near exact EM method provides almost exact estimation of the second moments of both displacement and velocity states which finally converges to the solution of weak 3 Taylor method. 

\begin{figure}[H]
	\centering
	\includegraphics[width=0.8\textwidth]{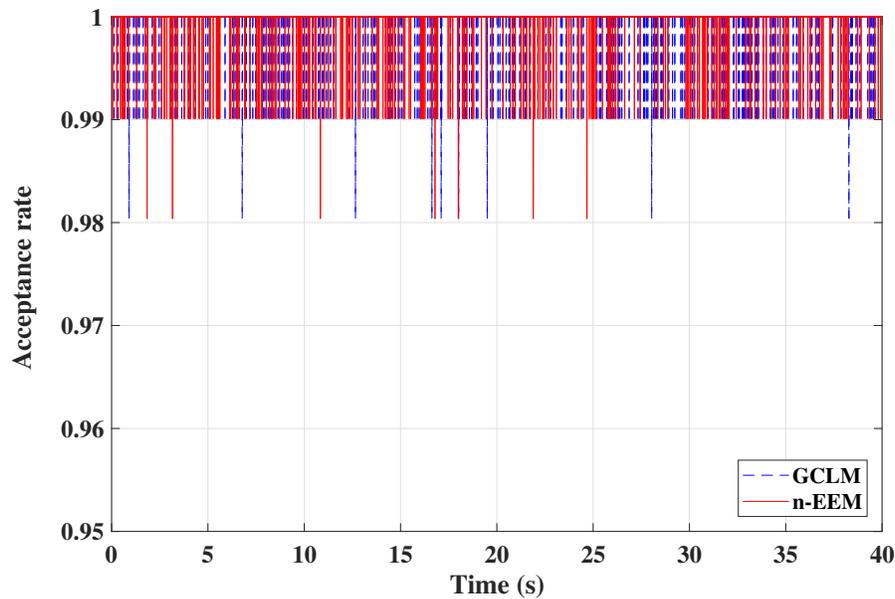}
	\caption{Acceptance ratio of the samples in rejection sampling for the 2-DOF oscillator}
	\label{fig:2dofv}
\end{figure}

But in cases of the existing GCLM scheme, the solutions only provides good approximation in case of second order displacement moments: $E\left[ {{X_1^2}} \right]$ and $E\left[ {{X_2^2}} \right]$, and shows a poor estimation in case of second order velocity moments. In Fig. \ref{fig:2dofv}, the acceptance ratio of the samples in rejection sampling for the proposed n-EEM scheme is compared with the GCLM method. Since, the solution of EM scheme without any correction lies in the proximity of the higher order solution, the acceptance ratios of the sample paths are almost 1.
\begin{figure}[H]
	\centering
	\includegraphics[width=0.8\textwidth]{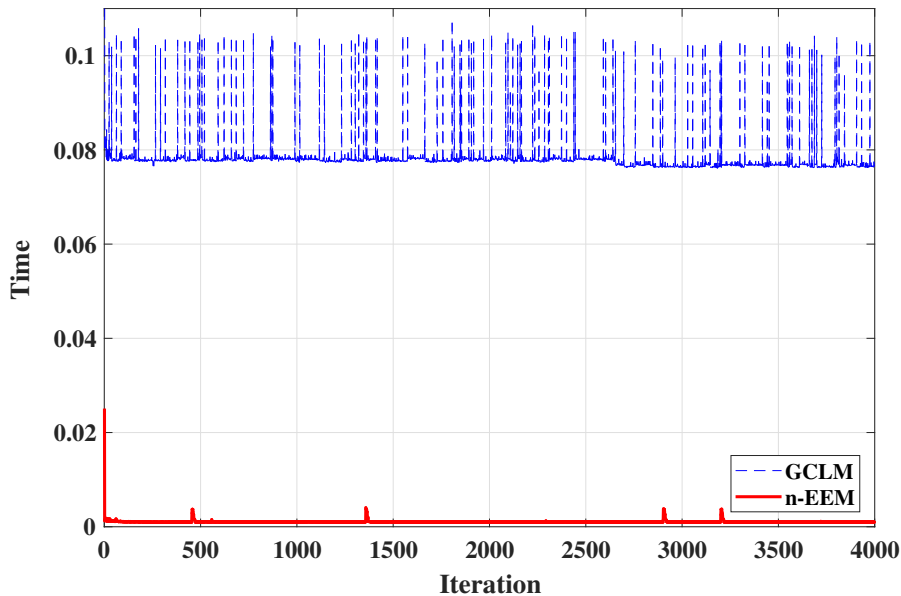}
	\caption{Amount of time involved in the rejection sampling to populate sample set at each time instant between the interval $t \in \left[0,40 \right]$s with $\Delta t = 0.01$}
	\label{fig:2doftime}
\end{figure}
Additionally it is observed that the time required for sample generation to populate the sample set in rejection sampling is almost identical for both the scheme. But, GCLM scheme exploits the local transversal scheme to generate the sample path for rejection sampling which is computationally expensive due to the need for computation of matrix exponential. The proposed n-EEM framework exploits a relatively efficient EM scheme to generate the sample paths for rejection sampling. Therefore, even for similar acceptance ratio, the proposed approach is computationally efficient as compared to the existing GCLM (see Fig. \ref{fig:2doftime}).

\section{Conclusions}
This work provides a new framework for the near exact simulation of the non-linear stochastic differential equations using Euler-Maruyama by formulating a change of measure based sampling strategy. In the simulation of non-linear systems, the errors mainly arise due to the limitation in the treatment of the non-linear functions rather than the linear function. In this work, the error in the non-linear approximation is eliminated from the solution by formulating a change of measure for the Brownian process through Girsanov's transform. The change of measure is computed in terms of a Radon-Nikodym derivative using a rejection sampling framework. The algorithm stipulates that the sample paths must satisfy the RN-derivative arsing due to the Girsanov's transformation. The efficacy of this proposed framework in solving nonlinear oscillators is studied using few non-linear stochastic systems and compared with stationary solutions wherever possible. In absence of a stationary solution, a comparison is made with higher order Weak 3.0 Taylor method. The method is further compared with a similar method namely GCLM. Since, the EM method is  effective and efficient than the local transversal scheme the acceptance ratio of sample paths in the proposed near exact Euler-Maruyama (n-EEM) scheme is naturally higher than the available Girsanov based LTL techniques. This provides good computational efficiency since very less time simulation is required to populate the sample path in rejection sampling. Further, whenever the available EM and Girsanov based LTL schemes fails to correctly approximate the higher order solution, the proposed n-EEM method provides almost near exact estimate for displacement and velocity states of a dynamical system. Thus it can be conjectured that in cases of a large variety of non-linear systems the proposed algorithm will provide an improvement over both classical EM and advanced Girsanov corrected techniques without involving high computational and derivational involvement like It\^{o}-Taylor Weak-3 scheme.

\appendix

\section{Stochastic exponential for RVP oscillator}\label{filter1}
The stochastic exponential: $\exp \left\{ {\int_{{t_{i - 1}}}^{{t_i}} {\gamma \left( {{{\bm X}_s},{{\bm X}_{i - 1}}} \right)d\tilde B\left( s \right)} } \right\}$ in Eq. (\ref{3.13}) for the RVP oscillator can be expanded noting that \\ $\gamma \left( {{{\bm X}_t},{{\bm X}_{i - 1}}} \right) =  - {\sigma ^{ - 1}}\left( \varepsilon \right)$, where, $\varepsilon  = {h_3}X_1^2(t){X_2}(t) + {h_3}X_2^3(t) - {h_3}X_{1,i - 1}^2{X_{2,i - 1}} - {h_3}X_{2,i - 1}^3$, using the It\^{o}-product rule as follows:
\begin{equation*}\label{4.20}
\begin{array}{ll}
\multicolumn{2}{l}{ \left. {\gamma \tilde B} \right|_{{t_{i - 1}}}^{{t_i}} =  - {\sigma ^{ - 1}}\left[ {\left( {{h_3}X_{1,i}^2{X_{2,i}} + {h_3}X_{2,i}^3} \right){{\tilde B}_i} - \left( {{h_3}X_{1,r}^2{X_{2,r}} + {h_3}X_{2,r}^3} \right){{\tilde B}_r}} \right.\left. { - \left( {{h_3}X_{1,r}^2{X_{2,r}} + {h_3}X_{2,r}^3} \right)\left( {{{\tilde B}_i} - {{\tilde B}_r}} \right)} \right] }\\
\multicolumn{ 2}{l}{ \left. {{{\tilde B}^2}\frac{{\partial \gamma }}{{\partial {X_2}}}\sigma } \right|_{{t_{i - 1}}}^{{t_i}} =  - \left[ {\tilde B_i^2\left( {{h_3}X_{1,i}^2 + 3{h_3}X_{2,i}^2} \right) - \tilde B_{i - 1}^2\left( {{h_3}X_{1,r}^2 + 3{h_3}X_{2,r}^2} \right)} \right] }\\
\int\limits_{{t_{i - 1}}}^{{t_i}} {\frac{{\partial \gamma }}{{\partial {X_2}}}\sigma ds}  = \int { - \left( {{h_3}X_1^2 + 3{h_3}X_2^2} \right)ds} , & \int\limits_{{t_{i - 1}}}^{{t_i}} {{{\tilde B}}\frac{{\partial \gamma }}{{\partial {X_1}}}{a_1}ds}  = \int {\tilde B\left[ { - {\sigma ^{ - 1}}\left( {2{h_3}{X_1}{X_2}} \right)} \right]{X_2}ds} \\
\int\limits_{{t_{i - 1}}}^{{t_i}} {{{\tilde B}}\frac{{\partial \gamma }}{{\partial {X_2}}}\tilde gds}  = \int {\tilde B\left[ { - {\sigma ^{ - 1}}\left( {{h_3}X_1^2 + 3{h_3}X_2^2} \right)} \right]\tilde gds}, & \frac{1}{2}\int\limits_{{t_{i - 1}}}^{{t_i}} {\tilde B\left( {\frac{{{\partial ^2}\gamma }}{{\partial X_2^2}}\sigma } \right)\sigma ds}  = \int {\frac{1}{2}\tilde B\left( { - 6{h_3}{X_2}} \right)\sigma ds} \\
\int\limits_{{t_{i - 1}}}^{{t_i}} {\left( {{{\tilde B}^2}\frac{{{\partial ^2}\gamma }}{{\partial {X_1}\partial {X_2}}}\sigma } \right){a_1}ds}  = \int {{{\tilde B}^2}\left( { - 2{h_3}{X_1}} \right){X_2}ds}, & \int\limits_{{t_{i - 1}}}^{{t_i}} {{{\tilde B}^2}\left( {\frac{{{\partial ^2}\gamma }}{{\partial X_2^2}}\sigma } \right)\tilde gds}  = \int {{{\tilde B}^2}\left( { - 6{h_3}{X_2}} \right)\tilde gds} \\
\int\limits_{{t_{i - 1}}}^{{t_i}} {\tilde B\left( {\frac{{{\partial ^2}\gamma }}{{\partial X_2^2}}\sigma } \right)\sigma ds}  = \int {\tilde B\left( { - 6{h_3}{X_2}} \right)\sigma ds} , & \int\limits_{{t_{i - 1}}}^{{t_i}} {{{\tilde B}^2}\left( {\frac{{{\partial ^2}\gamma }}{{\partial X_2^2}}\sigma } \right)\sigma d{B_s}}  = 0
\end{array}
\end{equation*}

\section{Stochastic exponential for DVP oscillator}\label{filter2}
Similar to the previous case, the stochastic components in the Radon-Nikodym derivative in Eq. (\ref{3.13}) for this oscillator can be alternatively evaluated by noting, $\gamma \left( {{{\bm X}_t},{{\bm X}_{i - 1}}} \right) =  - {\left( {\frac{{\rho {X_1}(t)}}{m}} \right)^{ - 1}}\left( {\frac{\alpha }{m}X_1^3(t) - \frac{\alpha }{m}X_{1,i - 1}^3} \right)$ and $\varepsilon  = \frac{\alpha }{m}X_1^3(t) - \frac{\alpha }{m}X_{1,i - 1}^3$ as follows:
\begin{equation*}\label{4.33}
\begin{array}{ll}
\multicolumn{2}{l}{ \left. {\gamma \tilde B} \right|_{{t_{i - 1}}}^{{t_i}} =  - {\left( {\frac{{\rho {X_1}(t)}}{m}} \right)^{ - 1}}\left( {\frac{\alpha }{m}X_{1,i}^3{{\tilde B}_i} - \frac{\alpha }{m}X_{1,r}^3{{\tilde B}_r} - \frac{\alpha }{m}X_{1,i - 1}^3\left( {{{\tilde B}_i} - {{\tilde B}_r}} \right)} \right)}\\
\multicolumn{2}{l}{ \left. {{{\tilde B}^2}\frac{{\partial \gamma }}{{\partial {X_2}}}\sigma } \right|_{{t_{i - 1}}}^{{t_i}} = 0}\\
\int\limits_{{t_{i - 1}}}^{{t_i}} {\frac{{\partial \gamma }}{{\partial {X_2}}}\sigma ds}  = 0, & \int\limits_{{t_{i - 1}}}^{{t_i}} {\tilde B\frac{{\partial \gamma }}{{\partial {X_1}}}{a_1}ds}  = \int {\tilde B\left( { - \frac{{3\alpha }}{m}{{\left( {\frac{{\rho {X_1}(t)}}{m}} \right)}^{ - 1}}X_1^2} \right){X_2}ds} \\
\int\limits_{{t_{i - 1}}}^{{t_i}} {\tilde B\frac{{\partial \gamma }}{{\partial {X_2}}}\tilde gds}  = 0, & \frac{1}{2}\int\limits_{{t_{i - 1}}}^{{t_i}} {\tilde B\left( {\frac{{{\partial ^2}\gamma }}{{\partial X_2^2}}\sigma } \right)\sigma ds}  = 0\\
\int\limits_{{t_{i - 1}}}^{{t_i}} {\left( {{{\tilde B}^2}\frac{{{\partial ^2}\gamma }}{{\partial {X_1}\partial {X_2}}}\sigma } \right){a_1}ds}  = 0, & \int\limits_{{t_{i - 1}}}^{{t_i}} {{{\tilde B}^2}\left( {\frac{{{\partial ^2}\gamma }}{{\partial X_2^2}}\sigma } \right)\tilde gds}  = 0\\
\int\limits_{{t_{i - 1}}}^{{t_i}} {\tilde B\left( {\frac{{{\partial ^2}\gamma }}{{\partial X_2^2}}\sigma } \right)\sigma ds}  = 0, & \int\limits_{{t_{i - 1}}}^{{t_i}} {{{\tilde B}^2}\left( {\frac{{{\partial ^2}\gamma }}{{\partial X_2^2}}\sigma } \right)\sigma d{B_s}}  = 0
\end{array}
\end{equation*}

\section{Stochastic exponential for 2-DOF oscillator}\label{filter3}
In this case, it is a straightforward extension of previous calculations to two-dimension. Noting that, \\ ${\gamma _1}\left( {{{\bm Y}_t},{{\bm Y}_{i - 1}}} \right) =  - \sigma _1^{ - 1}\alpha \left( {Y_1^2(t){Y_2}(t) - Y_1^2({t_{i - 1}}){Y_2}({t_{i - 1}})} \right)$ and ${\gamma _2}\left( {{{\bm Y}_t},{{\bm Y}_{i - 1}}} \right) =  - \sigma _2^{ - 1}\beta \left( {Y_3^3(t) - Y_3^3({t_{i - 1}})} \right)$, the followings are obtained:
\begin{equation*}
\begin{array}{ll}
\multicolumn{2}{l}{ \left. {{\gamma _1}{{\tilde B}_1}} \right|_{{t_{i - 1}}}^{{t_i}} =  - \sigma _1^{ - 1}\alpha \left( {Y_{1,i}^2{Y_{2,i}}{{\tilde B}_{1,i}} - Y_{1,r}^2{Y_{2,r}}{{\tilde B}_{1,r}} - Y_{1,r}^2{Y_{2,r}}\left( {{{\tilde B}_{1,i}} - {{\tilde B}_{1,r}}} \right)} \right) }\\
\multicolumn{2}{l}{ \left. {{\gamma _2}{{\tilde B}_2}} \right|_{{t_{i - 1}}}^{{t_i}} =  - \sigma _2^{ - 1}\beta \left( {Y_{3,i}^3 - Y_{3,r}^3 - Y_{3,r}^3\left( {{{\tilde B}_{1,i}} - {{\tilde B}_{1,r}}} \right)} \right) }\\
\left. {\tilde B_1^2\frac{{\partial {\gamma _1}}}{{\partial {Y_2}}}{\sigma _1}} \right|_{{t_{i - 1}}}^{{t_i}} =  - \alpha \left[ {\tilde B_{1,i}^2Y_{1,i}^2 - \tilde B_{1,r}^2Y_{1,r}^2} \right], & \left. {\tilde B_2^2\frac{{\partial {\gamma _2}}}{{\partial {Y_4}}}{\sigma _2}} \right|_{{t_{i - 1}}}^{{t_i}} = 0\\
\int\limits_{{t_{i - 1}}}^{{t_i}} {\frac{{\partial {\gamma _1}}}{{\partial {Y_2}}}{\sigma _1}ds}  = \int\limits_{{t_{i - 1}}}^{{t_i}} {\left( { - \alpha Y_1^2(t)} \right)ds} , & \int\limits_{{t_{i - 1}}}^{{t_i}} {\frac{{\partial {\gamma _2}}}{{\partial {Y_4}}}{\sigma _2}ds}  = 0,\\
\int\limits_{{t_{i - 1}}}^{{t_i}} {{{\tilde B}_1}\frac{{\partial {\gamma _1}}}{{\partial {Y_1}}}{g_1}ds}  = \int {{{\tilde B}_1}\left( { - 2\sigma _1^{ - 1}\alpha {Y_1}Y_2^2} \right)ds} , & \int\limits_{{t_{i - 1}}}^{{t_i}} {{{\tilde B}_2}\frac{{\partial {\gamma _2}}}{{\partial {Y_3}}}{g_3}ds}  = \int {{{\tilde B}_2}\left( { - 3\sigma _2^{ - 1}\beta Y_3^2{Y_4}} \right)ds} \\
\int\limits_{{t_{i - 1}}}^{{t_i}} {{{\tilde B}_1}\frac{{\partial {\gamma _1}}}{{\partial {Y_2}}}{{\tilde g}_1}ds}  = \int\limits_{{t_{i - 1}}}^{{t_i}} {{{\tilde B}_1}\left( { - \sigma _1^{ - 1}\alpha Y_1^2\tilde g} \right)ds} , & \int\limits_{{t_{i - 1}}}^{{t_i}} {{{\tilde B}_2}\frac{{\partial {\gamma _2}}}{{\partial {Y_4}}}{{\tilde g}_2}ds}  = 0,\\
\frac{1}{2}\int\limits_{{t_{i - 1}}}^{{t_i}} {{{\tilde B}_1}\left( {\frac{{{\partial ^2}{\gamma _1}}}{{\partial Y_2^2}}{\sigma _1}} \right){\sigma _1}ds}  = 0, & \frac{1}{2}\int\limits_{{t_{i - 1}}}^{{t_i}} {{{\tilde B}_2}\left( {\frac{{{\partial ^2}{\gamma _2}}}{{\partial Y_4^2}}{\sigma _2}} \right){\sigma _2}ds}  = 0\\
\int\limits_{{t_{i - 1}}}^{{t_i}} {\left( {\tilde B_1^2\frac{{{\partial ^2}{\gamma _1}}}{{\partial {Y_1}\partial {Y_2}}}{\sigma _1}} \right){g_1}ds}  = \int\limits_{{t_{i - 1}}}^{{t_i}} {\tilde B_1^2\left( { - 2\alpha {Y_1}} \right){Y_2}ds}, & \int\limits_{{t_{i - 1}}}^{{t_i}} {\left( {\tilde B_2^2\frac{{{\partial ^2}{\gamma _2}}}{{\partial {Y_3}\partial {Y_4}}}{\sigma _2}} \right){g_3}ds}  = 0,\\
\int\limits_{{t_{i - 1}}}^{{t_i}} {\tilde B_1^2\left( {\frac{{{\partial ^2}{\gamma _1}}}{{\partial Y_2^2}}{\sigma _1}} \right){{\tilde g}_1}ds}  = 0, & \int\limits_{{t_{i - 1}}}^{{t_i}} {\tilde B_2^2\left( {\frac{{{\partial ^2}{\gamma _2}}}{{\partial Y_4^2}}{\sigma _2}} \right){{\tilde g}_2}ds}  = 0\\
\int\limits_{{t_{i - 1}}}^{{t_i}} {{{\tilde B}_1}\left( {\frac{{{\partial ^2}{\gamma _1}}}{{\partial Y_2^2}}{\sigma _1}} \right){\sigma _1}ds}  = 0, & \int\limits_{{t_{i - 1}}}^{{t_i}} {{{\tilde B}_2}\left( {\frac{{{\partial ^2}{\gamma _2}}}{{\partial Y_4^2}}{\sigma _2}} \right){\sigma _2}ds}  = 0,\\
\int\limits_{{t_{i - 1}}}^{{t_i}} {\tilde B_1^2\left( {\frac{{{\partial ^2}{\gamma _1}}}{{\partial Y_2^2}}{\sigma _1}} \right){\sigma _1}d{B_1}} (s) = 0, & \int\limits_{{t_{i - 1}}}^{{t_i}} {\tilde B_2^2\left( {\frac{{{\partial ^2}{\gamma _2}}}{{\partial Y_4^2}}{\sigma _2}} \right){\sigma _2}d{B_2}(s)}  = 0
\end{array}
\end{equation*}
Identifying the errors: ${\varepsilon _1}(t) = Y_1^2(t){Y_2}(t) - Y_1^2({t_{i - 1}}){Y_2}({t_{i - 1}})$ and ${\varepsilon _2}(t) = Y_3^3(t) - Y_3^3({t_{i - 1}})$ it can be found that: ${\phi _1}(t) =  - \alpha Y_1^2(t) - 2{{\tilde B}_1}(t)\sigma _1^{ - 1}\alpha {Y_1}(t)Y_2^2(t) - {{\tilde B}_1}(t)\sigma _1^{ - 1}\alpha Y_1^2(t)\tilde g - 2B_1^2(t)\alpha {Y_1}(t){Y_2}(t)$ and ${\phi _2}(t) =  - 3\sigma _2^{ - 1}\beta Y_3^2(t){Y_4}(t)$. Then, one can find the $\phi (t)$ though the relations: $\varepsilon (t) = \sum\nolimits_{k = 1}^n {{\varepsilon _k}(t)}$ and $\phi (t) = \sum\nolimits_{k = 1}^n {{\phi _k}(t)}$ as:
\begin{equation*}
\phi (t) =  - \alpha Y_1^2(t) - 2{{\tilde B}_1}(t)\sigma _1^{ - 1}\alpha {Y_1}(t)Y_2^2(t) - {{\tilde B}_1}(t)\sigma _1^{ - 1}\alpha Y_1^2(t)\tilde g - 2B_1^2(t)\alpha {Y_1}(t){Y_2}(t) - 3\sigma _2^{ - 1}\beta Y_3^2(t){Y_4}(t)
\end{equation*}

\textbf{Acknowledgements: } SC acknowledges the financial support received from IIT Delhi in form of seed grant.


\section*{Declarations}

\subsection*{Funding} The corresponding author received funding from IIT Delhi in form of seed grant.

\subsection*{Conflicts of interest} The authors declare that they have no conflict of interest.

\subsection*{Availability of data and material} The datasets generated during and/or analysed during the current study are available from the corresponding author on reasonable request.

\subsection*{Code availability} The MATLAB codes written for this work are available from the corresponding author on reasonable request.


\end{document}